\input amstex
\documentstyle{amsppt} \magnification=1000
\hcorrection{0mm} \vcorrection{0mm} \hsize=13.5cm \vsize=21.8cm
\pageno=1 \NoBlackBoxes \nopagenumbers \nologo \headline={ \rightheadline }

\def\bproc#1{\parindent=0mm \bf #1.\hskip3mm\it}
\def\eproc{\rm \vskip1mm \parindent=4mm}
\def\proof#1{ \parindent=0mm {{\it Proof}\, #1:\hskip4mm}\parskip=1mm}
\def\qed{\newline \rightline{$\square$} \vskip3mm \parskip=5mm \parindent=4mm}
\def\bitem{\vskip2mm \leftskip8mm \parindent-8mm \parskip=2mm}
\def\eitem{\vskip-1mm \leftskip0mm \parindent4mm \parskip=5mm}
\def\rmk#1{{\noindent \bf #1.\hskip2mm}}
\def\prg#1#2{\vskip8mm \noindent {\S #1 \quad \bf #2.} \vskip1mm}

\def\trm#1{ {\text{\rm{#1}}} }
\def\ra{\rightarrow} \def\lra{\longrightarrow} \def\Lra{\Longrightarrow} \def\ms{\mapsto}
\def\Pic{{\operatorname{Pic}}} \def\modulo{{\operatorname{modulo}}}
\def\dim{{\operatorname{dim}\pt}}     
\def\mult{{\operatorname{mult}}}

\def\mmm#1#2#3{{ \pmatrix #1 &  & 0 \\  & #2 &  \\ 0 &  & #3 \endpmatrix }}

\def\blr#1#2{ { \buildrel{#2} \over {#1} } }


\def\mls#1{\operatorname{mult}({#1})_{\star}} \def\mus#1{\operatorname{mult}({#1})^{\star}}
\def\smr{{\trm{sum}^{\star}}}

\def\ptr{\pt \bigstar \pt}
\def\nts{\pt {}_{{}_\bullet} \pt}

\def\sqb#1#2{\pt\left|{\matrix{#1}\\ {#2}\endmatrix} \right|\pt}

 \def\Z{\Bbb Z} \def\N{\Bbb N} \def\Q{\Bbb Q} \def\P{\Bbb P}
\def\CC{{\Cal C}} \def\FF{{\Cal F}}  \def\PP{{\Cal P}}  \def\JJ{{\Cal J}}
\def\SS{{\Cal S}} \def\AA{{\Cal A}} \def\RR{{\Cal R}}
\def\pt{\hskip1pt}

\hyphenation{pola-rized} \hyphenation{Kummer} \hyphenation{theorem} \hyphenation{variety} \hyphenation{bundles}
\hyphenation{introduce} \hyphenation{property}

\def\bed{1} \def\fut{2} \def\tgf{3} \def\ntt{4} \def\fte{5} \def\okt{6} \def\ppz{7} \def\csp{8}
\def\ppt{9} \def\pzn{10} \def\pzl{11} \def\pzm{12} \def\pzc{13} \def\pzr{14} \def\pzd{15} \def\pzv{16}
\def\rbp{17} \def\mpt{18} \def\mpg{19} \def\mpb{20} \def\cij{21} \def\rij{22} \def\rijb{23}
\def\cvn{24} \def\chij{25} \def\ami{26} \def\als{27} \def\alsb{28} \def\alsc{29}  \def\rbe{30}
\def\lma{31} \def\lmb{32} \def\lmc{33}

\document
\parskip=5mm \parindent=4mm

 at 8truept
\font\bigfont=cmr12 at 18truept


\ \vskip15mm

\baselineskip 24pt
\centerline{\bigfont Tautological Cycles on Jacobian Varieties}

\baselineskip 14pt \vskip6mm

\bf \noindent
Giambattista Marini \newline
University of Rome \lq\lq Tor Vergata" \newline
00133 Roma (Italy)
\rm \vskip6mm

\prg{1}{Introduction}

Let $ \, \CC \, $ be a complex curve of genus $ \, g \, $ and let $ \, J(\CC) \, $ denote its Jacobian.
We consider the group of rational cycles modulo algebraic equivalence
$$
\AA_{\bullet}\big(J(\CC)\big)_{\Q}
$$
and the so called \lq\lq tautological ring" $ \, \RR(\CC) \, , \ $ that is, the subgroup of $ \, \AA_{\bullet}\big(J(\CC)\big)_{\Q} \, $
containing $ \, \CC \, $ and stable with respect to the Fourier transform, the intersection product, the Pontryagin product,
pull-backs and push-forwards of multiplication maps by integers.
The tautological ring has a mysterious algebraic structure, which has not been completely understood so far.
Only recently Beauville proved that $ \, \RR(\CC) \, $ is finite dimensional as a $ \, \Q $-vector space [Be3].
He proved that the $ \, W^d $'s and their intersections generate
$ \, \RR(\CC) \, $ as a vector space (this gives a rough bound for its dimension).
The difficulty of understanding the structure of $ \, \RR(\CC) \, $ comes from the fact that
for a Jacobian, the cycles occurring in nature, such as $ \, \CC \, $ itself or the $ \, W^d $'s
(i.e. the varieties parameterizing
effective $ d $-degrees line bundles on $ \, \CC \, , \ $ or, equivalently, the products $ \, \CC^{\ptr d} / d! \, ) \ $
are not pure with respect to Beauville's graduation [Be2].
From the point of view of Fourier duality, the products of components of $ \, \CC \, , \ $ though they are
pure cycles with respect to both graduations, do not have a nice behavior under the intersection with the Theta divisor.
This creates trouble when trying to describe the tautological ring of the curve $ \, \CC \, : \ $
the matrices representing the intersection product and the Pontryagin product look completely chaotic if one works with
such cycles, and it is difficult to even determine an exhaustive set of relations among them.

In this paper we study the algebraic structure of $ \, \RR(\CC) \, , \ $
giving a detailed description of all the possibilities that may occur
(for $ g \le 8 ): \ $ we construct convenient basis and we determine the
matrices representing the Fourier transform and both intersection and Pontryagin products explicitly.
In particular, we estimate the dimension of $ \, \RR(\CC) \, $. \
We also define algebraic models for $ \, \RR(\CC) \, . \ $
We obtain our descriptions as a consequence of a more general result, that can be applied also to higher genus cases.

In the next section, we recall some basic results and fix the notation.
Section 3 is devoted to technical results concerning the Fourier transform of a pure cycle whose intersection with the
Theta divisor is trivial.
In section 4, we construct cycles for the tautological $ \, \RR(\CC) \, $ ring of a Jacobian, which have a nice behavior in some sense.
Section 5 is specular to section 4: we introduce abstract algebraic models for $ \, \RR(\CC) \, . \ $
In section 6, we apply our machinery to give an exhaustive description of $ \, \RR(\CC) \, $ for
all the possibilities that may occur in the cases where $ \, g \le 8 \, . \ $
We also consider the cases where $ \, \CC \, $ (a curve of any genus) has at least a $ \, g_4^1 \, . $

\prg{2}{Preliminary lemmas \& definitions}

Throughout this section, $ \, (X,\Theta) \, $ denotes a principally polarized abelian variety (p.p.a.v.) of dimension $ \, g \, , \ $
$ \, \widehat X \, = \, \Pic^0(X) \, $ denotes the dual abelian variety of $ \, X \, $ and
$ \ \AA_{\bullet}(X)_{\Q} $ \ denotes the group of rational algebraic cycles modulo algebraic equivalence,
graded by dimension. We also set
$$
\AA^{g-d}(X)_{\Q} \quad = \quad \AA_d(X)_{\Q} \quad = \quad \trm{\lq\lq subgroup of $d$-dimensional cycles"}
$$
The group $ \ \AA_{\bullet}\big(X)_{\Q} \ $ has two ring structures:
the first one is given by the intersection product [Fu], the second one is given by the Pontryagin product \lq\lq$ \ptr $",
that we recall to be defined by $ \, \alpha \ptr \beta \, = \, \SS_{\star} (\alpha \times \beta) \, , \ $
where $ \, \alpha \times \beta \, $ is the natural product in $ \, X \times X \, $ and
$ \ \SS : \, X \times X \, \ra \, X \ $ denotes the sum map on $ \, X \, . \ $
The intersection product is homogeneous with respect to codimension and its unit is given by
the class of $ \, X \, $ itself. The Pontryagin product is homogeneous with respect to dimension and its unit is given by
$ \, \{o\} \, , \ $ where $ \, o \, $ denotes the origin of $ \, X \, $ and
$ \, \{o\} \, $ denotes the class of a point in $ \, \AA_0\big(X\big)_{\Q} \, . $
Besides the two products mentioned above, the group $ \ \AA_{\bullet}\big(X)_{\Q} \ $ has a remarkable automorphism, namely
the Fourier transform, and it has a second graduation defined in terms of pull-backs
$ \, \mus{m} \, . \ $ Here, $ \mult(m) $ denotes the endomorphism of $ X $ given by the multiplication map by the integer $ m . $
\ The $ \, d $-dimensional cycles of degree $ \, s \, $ are defined by
$$
\left[ \AA_d\big(X\big)_{\Q} \right]_s \hskip2mm := \hskip2mm
\left\{ W \in \AA_d\big(X\big)_{\Q} \bigg\vert \mus{m} W = m^{2g-2d-s} W , \ \forall \, m \in \Z \right\}
$$
(we use the same notation from [Be2] and shall refer to this graduation as Beauville's graduation).
It is worth noticing that this graduation can also be defined in terms of push-forwards $ \, \mls{m} \, : \ $
from [Be2], we know that
$ \, \mls{m} W \, = \, m^{2d+s} W \, , \ $ for all $ \, m \in \Z \, , \ $ if and only if
$ \, \mus{m} W \, = \, m^{2g-2d-s} W \, , \ $ for all $ \, m \in \Z \, . \ $
The effect on Chow groups of such pull-backs and push-forwards has been studied by
Beauville in [Be1] and [Be2], Dehninger and Murre in [DM] and K\"unnemann in [K].
Dimensional graduation and Beauville's graduation give rise to a direct sum decomposition which reduces to
$$
\AA^p\big(X\big)_{\Q} \quad = \quad
\bigoplus_{s \, = \, p-g}^{\max\{0, \, p-1\}} \ \left[ \AA^p\big(X\big)_{\Q} \right]_s \ ,
\tag{\bed.1}
$$
if we take into account the fact that we are working modulo algebraic equivalence
(see [Be2], main theorem and proposition $(4.a)). \ $ We shall use the fact that, in particular,
$$
\left[ \AA^p\big(X\big)_{\Q} \right]_s \ = \ 0 \ , \qquad s \ \ge \ \max\,\{p,\, 1\}
\tag{\bed.2}
$$
Clearly,
$ \ \AA^0\big(X\big)_{\Q} \, = \, \left[ \AA^0\big(X\big)_{\Q}\right]_0 \, = \, \Q X \, , \ $
\ $ \AA^1\big(X\big)_{\Q} \, = \, \left[ \AA^1\big(X\big)_{\Q}\right]_0 \, \cong \, \frac{\Pic(X)}{\Pic^0(X)} \, , $
\ $ \AA_0\big(X\big)_{\Q} \, = \, \left[ \AA_0\big(X\big)_{\Q}\right]_0 \, = \, \Q \{o\} \, . $

We now recall briefly the definition and the main features of the Fourier transform. Let
$$
\PP \quad = \quad \Theta \times X + X \times \Theta - \smr \Theta \quad \in \quad \AA^1(X\times \widehat X)_{\Q}
$$
denote the normalized Poincar\'e divisor. The Fourier transform is the map
$$
\eqalign{
\FF : \ \AA_{\bullet}(X)_{\Q} \ & \lra \ \AA_{\bullet}(\widehat X)_{\Q} \, ; \cr
W \hskip6mm & \ \ms \, \hskip3mm e^{\PP}(W)
}
$$
where $ \ e^{\PP} \ = \ X\!\!\times\!\!\widehat X + \PP + \PP^2/2 + \PP^3/3! + ... \ $ is the
exponential of $ \, \PP \, $ under the intersection product and $ \, e^{\PP} \, $ acts on $ \, W \, $ as a correspondence.
We identify $ \, \AA_{\bullet}(\widehat X)_{\Q} \, $ with $ \, \AA_{\bullet}(X)_{\Q} \, $ via the natural
identification of $ \, \widehat X \, $ with $ \, X \, , \ $ induced by the principal polarization $ \, \Theta \, . $
\ Under this identification, we consider the Fourier transform $ \, \FF \, $ as a transform on $ \, \AA_{\bullet}(X)_{\Q} \, . $
\ The main features of Fourier transform are the following (see [Muk], [Be1], [Be2] and [Be3]):

\bitem
$ (\fut.1) \qquad
\FF\big(\alpha \ptr \beta\big) \ = \ \FF(\alpha) \nts \FF(\beta) \ $ and
$ \ \FF\big(\alpha \nts \beta\big) \ = \ \big(-1\big)^g \FF(\alpha) \ptr \FF(\beta) \, ; $

$ (\fut.2) \qquad
\FF \circ \FF \ = \ \big(-1\big)^g \mus{-1} \, , \ $ in particular $ \, \FF \,$ is bijective;

$ (\fut.3) \qquad
\FF\,\left[ \AA_d\big(X\big)_{\Q} \right]_s \ = \ \left[ \AA^{d+s}\big(X\big)_{\Q} \right]_s \, ; $

$ (\fut.4) \qquad
\FF(W) \ = \ e^{\Theta} \nts \left[ \left(W^{-} \nts e^{\Theta} \right) \ptr e^{- \Theta}\right] \, , \ $ where
$ \, W^{-} := \mus{-1} W \, . $
\eitem

\noindent
Furthermore, both the intersection product and the Pontryagin product are homogeneous with respect to Beauville's graduation.

A remarkable 1-cycle which lives in 0-Beauville's degree (likewise the theta divisor)
is $ \ \Gamma \, = \, \frac{\Theta^{g-1}}{(g-1)!} \, . \ $
A few basic well known identities involving $ \, \Gamma \, $ are
$ \ \FF(\Theta) \, = \, (-1)^{g-1} \Gamma \, , \ \FF(\Gamma) \, = \, - \Theta \ $ and
$ \ \frac{\Theta^k}{k!} \, = \, \frac{\Gamma^{\ptr g-k}}{(g-k)!} \, . \ $
To our purpose, it is convenient to write explicitly two trivial consequences of the basic identities above.
From (\fut.1) and (\fut.2), it is a straightforward exercise to check that

\bitem
$ (\tgf.1) \qquad
\dsize \FF \left(\frac{\Gamma^{\ptr m}}{m\,!}\right) \quad = \quad \big(-1\big)^m \, \frac{\Gamma^{\ptr g-m}}{(g-m)!} \, ; $

$ (\tgf.2) \qquad
\dsize \left(\frac{\Gamma^{\ptr m}}{m\,!}\right) \nts \frac{\Theta^h}{h!} \quad = \quad
{g\!-\!m\!+\!h \choose h} \, \frac{\Gamma^{\ptr m-h}}{(m-h)\,!} \, . $
\eitem

\noindent
We may note that under the notation below (which we fix for the whole paper), the identities (\tgf.1) and (\tgf.2) hold
for all $ \, m \, $ and $ \, h \, $ in $ \, \Z \, . $

\rmk{Notation \ntt}
Let $ \, \alpha \, \in \, \AA_{\bullet}\big(X\big)_{\Q} \, $ be a cycle,  we set
$$
\alpha^k \ = \ \alpha^{\ptr k} \ = \  0 \ , \ \ \forall \ k < 0 \ ; \qquad
\alpha^0 \, = \, X \, , \ \ \alpha^{\ptr 0} \, = \, \{o\} \, , \ \ \forall \ \alpha \ne 0 \, .
$$
We also use the standard notation (defined for all $ \, n \in \Z \, , \ $ cf. [ACGH]):
$$
{n \choose h} \, = \, \frac{n \cdot (n\!-\!1)\cdot ... \cdot(n\!-\!h\!+\!1)}{h!} \, , \ h > 0 \ ; \qquad
{n \choose 0} \, = \, 1 \ ; \qquad
{n \choose h} \, = \, 0 \, , \ h < 0
$$
(and $ \ k\,! \ = \ 1 \, , \ \forall \ k \le 0 \, ). $


\vskip14mm

\prg{3}{On the Fourier transform}

This section deals with remarkable properties of pure cycles (dimensionally and with respect to Beauville's graduation)
having a trivial intersection with the theta divisor.
We show that such cycles define remarkable subspaces of the Chow ring:
for $ \ W \, \in \, \left[\AA_d\big(X\big)_{\Q}\right]_s $ \ satisfying $ \, W \nts \Theta \, = \, 0 \, $
and $ \, \Gamma \, $ as in \S 2, the $ \, \Q $-vector space with basis
$ \ \left\{ \, W \ptr \Gamma^{\ptr m} \right\}_{m=1}^{g-s-2d} $ \ is $ \, \FF $-stable and turns out to be of particular interest:
while dealing with Jacobians, such kind of cycles
will be essential to our description of the tautological ring $ \, \RR(\CC) \, . \ $

Throughout this section $ \, (X,\Theta) \, $ denotes a p.p.a.v. of dimension $ \, g \ge 1 \, $ and
$ \ W \, \in \, \AA_{\bullet}\big(X\big)_{\Q} $ \ denotes a cycle.
Writing formula (\fut.4) explicitly, we obtain
$$
\FF(W) \qquad = \qquad \sum_{t,\,h,\,j\ \ge\ 0} \quad
\left[ \left( W^{-} \nts \frac{\Theta^j}{j!}\right)\ptr \frac{\big(-\Theta\big)^h}{h!}\right] \nts \frac{\Theta^t}{t!}
\tag{\fte.1}
$$
Our next task is to simplify formula (\fte.1) in the case when $ \, W \, $ is pure,
(see formula \okt.6 below). We start with a definition.

\rmk{Definition}
For $ \ W \, \in \, \left[ \AA_d\big(X\big)_{\Q} \right]_s \ $ and $ \ k , \, t \ \in \ \Z \ $ we define
$$
\Omega_{k;\,t}(W)\quad = \quad \sum_{j=0}^d \ \big(-1\big)^{k+t+j} \
\left[ \left( W \nts \frac{\Theta^j}{j!}\right)\ptr \frac{\Theta^{s+2d-j-t-k}}{(s+2d-j-t-k)!}\right] \nts \frac{\Theta^t}{t!} \ .
$$

We may note that for $ \ W \, \in \, \left[ \AA_d\big(X\big)_{\Q} \right]_s \, , \ $
in consideration of the previous definition, of notation (\ntt) and also in consideration of
the identity $ \ W^- := \mus{-1} W = (-1)^s W \, , $ \ the equality (\fte.1) reduces to
$$
\FF(W) \qquad = \qquad \sum_{k , \, t \, \in \, \Z } \quad \Omega_{k;\,t}(W)
\tag{\fte.2}
$$
In the following lemma we prove some properties of $ \ \Omega_{k;\,t}(W) \, . $

\bproc{Lemma \okt}
Let $ \, (X,\Theta) \, $ be a p.p.a.v. of dimension $ \, g \ge 1 \, , \ $ let
$ \ W \, \in \, \left[ \AA_d\big(X\big)_{\Q} \right]_s \ $ and let $ \, \Omega_{k;\,t} = \Omega_{k;\,t}(W) \, . $ \ Then

\bitem
$(\okt.1)$ \hskip14mm $ \dim \Omega_{k;\,t} \quad = \quad g-d-s+k \, ; $

$(\okt.2)$ \hskip14mm $ \dsize\sum_{t\ge 0} \ \Omega_{k;\,t} \quad = \quad 0 \, , \quad \forall \ k \ne 0 \, ; $

$(\okt.3)$ \hskip14mm $ \Omega_{k;\,t} \quad = \quad \left(\Omega_{k+t;\,0}\right) \nts \dsize\frac{\Theta^t}{t!} \, ,
\quad \forall \ t \in \Z $ \hfill {\rm (both sides vanish for } $ t < 0);$

$(\okt.4)$ \hskip14mm $ \Omega_{k;\,t} \quad = \quad 0 \, , \quad \forall \ k+t \ge 1 \, ; $

$(\okt.5)$ \hskip14mm $ \Omega_{-m;\,0} \quad = \quad \left(\Omega_{0;\,0}\right) \nts \dsize\frac{(-\Theta)^m}{m!} \, , \quad \forall \ m \ge 0 \, ; $

$(\okt.6)$ \hskip14mm $ \FF(W) \quad = \quad \Omega_{0;\,0} \quad = \quad
\dsize\sum_{j=0}^d \ (-1)^j \left( W \nts \frac{ \Theta^j}{j!} \right) \ptr \frac{ \Theta^{s+2d-j}}{(s+2d-j)!} \ . $
\eitem
\eproc

Note that, in particular, for a 1-cycle $ \, Y \, \in \, \left[ \AA_1\big(X\big)_{\Q} \right]_s \, $ with $ \, s \ge 1 \, , \ $
the intersection $ \ Y \nts \Theta \ $ is trivial and formula $ \, (\okt.6) \, $ reduces to the identity
$ \ \FF(Y) \, = \, Y \ptr \frac{\Theta^{s+2}}{(s+2)!} \, . $

\proof{}
Equality (\okt.1) is obvious.
In view of (\fte.2) and (\okt.1), equality (\okt.2) follows from the fact that $ \ \FF(W) \ $ is pure of dimension $ \, g-d-s \, $
(see \fut.3). Furthermore, equality (\okt.3) is straightforward.

We now prove equality (\okt.4).
First, in view of (\okt.3), it suffices to prove that $ \ \Omega_{k;\,0} \, = \, 0 \, , \ \forall \ k \, \ge \, 1 \, . \ $
Furthermore, as $ \ \Omega_{k;\,0} \, \in \, \left[ \AA^{s+d-k}\big(X\big)_{\Q} \right]_s \, , \ $ the cases where
$ \ k > d \ $ follow from (\bed.2) for $ \ s > 0 \, , \ $ and they follow from trivial codimension reasons for $ \ s \le 0 \, . $

For $ \ 1 \le k \le d \ $ we proceed by descending induction on $ \, k : \ $
fix $ \, k \, $ in the range considered and assume
we are done for all $ \, k' > k \, . \ $ We have
$$
0 \quad
= \quad \sum_{t\ge 0} \ \Omega_{k;\,t} \quad
= \quad \sum_{t\ge 0} \ \left(\Omega_{k+t;\,0}\right) \nts \frac{\Theta^t}{t!} \quad
= \quad \Omega_{k;\,0}
\tag{\okt.7}
$$
where the first equality follows by (\okt.2), the second one follows by (\okt.3)
and the third equality follows by the inductive hypothesis
(the only term from the summation that survive is the one where $ \, t = 0). $

We now prove equality (\okt.5). For $ \, m = 0 \, $ (\okt.5) is trivial. Let $ \, m \ge 1 \, . \ $
First of all, we note the following
$$
0 \quad = \quad \sum_{t\ge 0} \ \Omega_{-m;\,t} \quad = \quad
\sum_{t=0}^m \ \Omega_{-m;\,t} \quad = \quad
\sum_{t=0}^m \ \left(\Omega_{-m+t;\,0}\right) \nts \frac{\Theta^t}{t!} \ ,
\tag{\okt.8}
$$
where the first equality follows from (\okt.2), the second identity follows from (\okt.4)
and the third one follows from (\okt.3). \
We conclude by induction on $ \, m \, . \ $
Assume we are done for all $ \, m' \, $ satisfying $ \, 0 \le m' < m \, . \ $ We have
$$
\eqalign{
\Omega_{-m;\,0} \quad = \quad - \ \sum_{t=1}^m \ \big[\Omega_{-m+t;\,0}\big] \nts \frac{\Theta^t}{t!} \quad
& = \quad - \ \sum_{t=1}^m \ \left[\big(\Omega_{0;\,0}\big) \nts \frac{\big(-\Theta\big)^{m-t}}{(m-t)!} \right] \nts \frac{\Theta^t}{t!} \cr
& = \quad - \ \big(\Omega_{0;\,0}\big) \nts \frac{\big(-\Theta\big)^m}{m!} \
\left[\sum_{t=1}^m \big(-1\big)^t \frac{m!}{(m-t)! \, t!} \right] \cr
& = \quad \big(\Omega_{0;\,0}\big) \nts \frac{\big(-\Theta\big)^m}{m!}
}
$$
where the first identity follows from (\okt.8), the second one follows by the inductive hypothesis, the third one and the forth one
are straightforward.

We now prove equality (\okt.6).
We have $ \ \FF(W) \, = \, \sum_{t\ge0} \Omega_{0;\,t} \ = \ \Omega_{0;\,0} \, , \ $
where the first equality follows from (\fte.2) and (\okt.2), the second equality follows from (\okt.4).
At last, $ \ \Omega_{0;\,0} \ $ equals the sum at the right hand side of (\okt.6).
\qed

\bproc{Proposition \ppz}
Let $ \, (X,\Theta) \, $ be a p.p.a.v. of dimension $ \, g \ge 1 \, . \ $
Consider $ \ W \, \in \, \left[\AA_d\big(X\big)_{\Q} \right]_s \ $ and assume
$ \ W \ne 0 \, , \ W \nts \Theta = 0 \, . \ $ Then,

\bitem
$ (\ppz.1) \qquad \dsize 0 \, \le \, s+2d \, \le \, g \ ; $

$ (\ppz.2) \qquad \dsize W \ptr \big( \Theta^{g-m} \big) \quad \ne \quad 0 \qquad
\Longleftrightarrow \qquad 0 \ \le \ m \ \le \ g-s-2d \ . $
\eitem

\noindent
Furthermore, for all $ \ h , \, m \, \in \, \Z \ $ the following holds

\bitem
$ (\ppz.3) \qquad \dsize \FF \left( W \ptr \frac{\ \Theta^{g-m} }{(g-m)!} \right) \quad
= \quad (-1)^m \ W \ptr \frac{\Theta^{s+2d+m}}{(s+2d+m)!} \ ; $

$ (\ppz.4) \qquad \dsize \left[W \ptr \frac{\ \Theta^{g-m} }{(g-m)!}\right] \nts \frac{\Theta^h}{h!} \
= \ {g\!-\!s\!-\!2d\!-\!m\!+\!h \choose h} \ W \ptr \frac{\ \Theta^{g-m+h} }{(g-m+h)!} \ ; $
\eitem
\eproc

\proof{}
First, we observe that as $ \ W \nts \Theta \, = \, 0 \, , \ $ the cycle $ \, \Omega_{k; \, t}(W) \, $ reduces to
$$
\Omega_{k; \, t}(W) \quad = \quad \big(-1\big)^{k+t} \ \left[ W \ptr \frac{\Theta^{s+2d-k-t}}{(s+2d-k-t)!} \right] \nts \frac{\Theta^t}{t!}
$$
Then, applying (\okt.6) we obtain
$ \ \FF(W) \, = \, \Omega_{0;\, 0}(W) \, = \, W \ptr \frac{\Theta^{s+2d}}{(s+2d)!} \, . \ $
In particular, (\ppz.1) holds \ (for otherwise the Fourier transform of $ \, W \, $ would be trivial,
contradicting the hypothesis $ \, W \ne 0). \ $

{\it Step 1. \ } We prove (\ppz.4) in the special case $ \, g-m \, = \, s+2d \, . \ $
As (\ppz.4) is trivial for $ \, h < 0 \, , \ $ this reduces to prove the following
$$
\left[W \ptr \frac{\Theta^{s+2d}}{(s+2d)!}\right] \nts \frac{\Theta^h}{h!} \quad
= \quad W \ptr \frac{\Theta^{s+2d+h}}{(s+2d+h)!} \ , \qquad \forall \ h \ge 0 \ .
\tag{\ppz.5}
$$
The case when $ \, h = 0 \, $ is trivial. Then we proceed by induction.
We fix $ \, h > 0 \, $ and assume we are done for all $ \, h' \, $ satisfying $ \, 0 \le h' < h \, . \ $
From (\okt.4), applying (\okt.2) with $ \, k = -h \, $ we obtain
$$
0 \quad = \quad \sum_{t=0}^h \ \Omega_{-h; \, t}(W) \quad =
\quad \sum_{t=0}^h \ (-1)^{t-h} \ \left[ W \ptr \frac{\Theta^{s+2d+h-t}}{(s+2d+h-t)!} \right] \nts \frac{\Theta^t}{t!} \ .
\tag{\ppz.6}
$$
For $ \, t > 0 \, $ we can apply the induction hypothesis to the term between brackets,
that is we can write
$ \, W \ptr \frac{\Theta^{s+2d+h-t}}{(s+2d+h-t)!} \, = \,
\left[ W \ptr \frac{\Theta^{s+2d}}{(s+2d)!} \right] \nts \frac{\Theta^{h-t}}{(h-t)!} \, . \ $
As a consequence, from (\ppz.6) we get
$$
\eqalign{
W \ptr \frac{\Theta^{s+2d+h}}{(s+2d+h)!} \quad & = \quad
- \sum_{t=1}^h \ (-1)^t \ \left[ W \ptr \frac{\Theta^{s+2d+h-t}}{(s+2d+h-t)!} \right] \nts \frac{\Theta^t}{t!} \cr
& = \quad
- \sum_{t=1}^h \ (-1)^t \ \left[ W \ptr \frac{\Theta^{s+2d}}{(s+2d)!} \right] \nts \frac{\Theta^{h-t}}{(h-t)!} \nts \frac{\Theta^t}{t!} \cr
& = \quad
- \sum_{t=1}^h \ (-1)^t \ {h \choose t} \ \left[ W \ptr \frac{\Theta^{s+2d}}{(s+2d)!} \right] \nts \frac{\Theta^{h}}{h!} \cr
& = \quad
\left[ W \ptr \frac{\Theta^{s+2d}}{(s+2d)!} \right] \nts \frac{\Theta^{h}}{h!}
}
$$
This proves (\ppz.5) and concludes the proof of step 1.

{\it Step 2. \ } We prove (\ppz.3).
As $ \ \FF(W) \, = \, W \ptr \frac{\Theta^{s+2d}}{(s+2d)!} \ $ and
$ \ \FF\left(\frac{\Theta^{g-m}}{(g-m)!}\right) \, = \, \frac{(-\Theta)^m}{m!} \, , $
\ applying (\fut.1) we obtain
$$
\eqalign{
\FF \left( W \ptr \frac{\ \Theta^{g-m} \ }{(g-m)!} \right) \quad
& = \quad (-1)^m \ \left[W \ptr \frac{\Theta^{s+2d}}{(s+2d)!}\right] \nts \frac{\Theta^m}{m!} \cr \cr
& = \quad (-1)^m \ W \ptr \frac{\Theta^{s+2d+m}}{(s+2d+m)!}
}
\tag{\ppz.7}
$$
where the second equality holds a priori only for $ \, m \ge 0 \ $ (apply \ppz.5 with $ \, h = m ). \ $
This proves (\ppz.3) for $ \, m \ge 0 \, . \ $

On the other hand, for $ \, m > g-s-2d \, $ (which is non-negative by \ppz.1) the right-hand-side of (\ppz.3) vanishes
by trivial dimensional argument.
Thus, the left-hand-side of (\ppz.3) must vanish as well, and, since the Fourier transform is injective, the following holds
$$
W \ptr \Theta^{g-m} \quad = \quad 0 \ , \qquad m \ > \ g-s-2d \, .
\tag{\ppz.8}
$$
This vanishing is exactly the vanishing of the right hand side of (\ppz.3) for negative values of $ \, m \, . \ $
The left hand side of (\ppz.3) vanishes for $ \, m < 0 \, $ as well (for trivial dimensional argument).
This concludes the proof of (\ppz.3).

{\it Step 3. \ } We prove (\ppz.2). By (\ppz.8),
we only have to prove that $ \ W \ptr \Theta^{g-m} \ $ is non-trivial for $ \, m \, $ in the range
$ \ 0 \le m \le g-s-2d \, . $ \ Let $ \, m \, $ be in this range and let
$ \, \mu = \frac{m!\ (g-s-2d-m)!}{(g-m)!\ (g-s-2d)!\ (m+s+2d)!} \, \Theta^{m+s+2d} \, . \ $ Since
$$
W \ptr \Theta^{g-m} \ptr \mu \quad = \quad W \ptr \frac{\Theta^{s+2d}}{(s+2d)!} \quad = \quad \FF(W) \quad \ne \quad 0 \ ,
$$
then $ \ W \ptr \Theta^{g-m} \ $ cannot vanish.

{\it Step 4. \ } We conclude the proof of (\ppz.4).
As already observed, (\ppz.4) is trivial for $ \, h \le 0 \, . \ $
So, let $ \, h > 0 \, . \ $
For $ \, g-s-2d-m < 0 \, , \ $ the left hand side of (\ppz.4) vanishes by (\ppz.2),
and the right hand side of (\ppz.4) vanishes simply because the coefficient
$ \, {g-s-2d-m+h \choose h} \, $ does vanish.
In the remaining case, namely $ \, h > 0 \, $ and $ \, g\!-\!s\!-\!2d\!-\!m \, \ge \, 0 \, , \ $ identity (\ppz.4) is
straightforward:
$$
\eqalign{
\left[W \ptr \frac{\ \Theta^{g-m} }{(g-m)!}\right] \nts \frac{\Theta^h}{h!} \quad
& = \quad \left[\left(W \ptr \frac{\Theta^{s+2d}}{(s+2d)!}\right) \nts \frac{\Theta^{g-s-2d-m}}{(g-s-2d-m)!}\right] \nts \frac{\Theta^h}{h!} \cr
& = \quad {g\!-\!s\!-\!2d\!-\!m\!+\!h \choose h} \left[\left(W \ptr \frac{\Theta^{s+2d}}{(s+2d)!}\right) \nts \frac{\Theta^{g-s-2d-m+h}}{(g-s-2d-m+h)!}\right] \cr
& = \quad {g\!-\!s\!-\!2d\!-\!m\!+\!h \choose h} \ W \ptr \frac{\ \Theta^{g-m+h} }{(g-m+h)!} \ ,
}
$$
where the second equality is trivial, the first and the third equality follow from (\ppz.5).
\qed

\newpage 

\prg{4}{The tautological ring $ \ \RR(\CC) \, $}

In this section $ \, \CC \, $ denotes a curve of genus $ \, g \, $ and
$$
\CC_{(s)} \quad \in \quad \left[ \AA_1\big(X\big)_{\Q} \right]_s
$$
denotes its $ \, s $-component from Beauville's decomposition.
We recall a result of Colombo-van Geemen [CG] stating that, assuming $ \, \CC \, $ is a $ \ m $-cover of $ \, \P^1 \, , \ $
$ \ \CC_{(s)} \ $ vanishes for $ s \ge m-1 \, , \ $ i.e.
$ \ \CC_{(s)} \, \ne \, 0 \ $ implies $ \ 0 \, \le \, s \, \le \, m-2 \, . \ $
Note that, as a consequence,
$$
\CC_{(s)} \ \ne \, 0 \qquad \Lra \qquad
0 \, \le \, s \, \le \, \frac{g-1}{2} \, . \
\tag{\csp}
$$
In fact, every curve of genus $ \, g \, $ is a $ \, m $-cover of $ \, \P^1 \, $ for some $ \, m \le \frac{g+3}{2} \ $ [ACGH].

We may note that $ \ \CC_{(0)} = \Gamma \ $ (cf. end of \S 2).
In fact, from (\fut.3), $ \ \FF(\CC_{(0)}) \, $ is the divisorial component of $ \, \FF(\CC) \, $ and, at the same time,
the divisorial component of $ \ \FF(\CC) $ \ is $ \ \, -\Theta \, = \, \FF(\Gamma) \, . \ $
Thus, identities (\tgf.1) and (\tgf.2) become

\bitem
$ (\ppt.1) \hskip5mm
\FF \left( \frac{{\CC_{(0)}}^{\ptr m}}{m!} \right) \ = \ (-1)^m \ \frac{{\CC_{(0)}}^{\ptr g-m}}{(g-m)!} \ , \quad m \in \Z \ ; $

$ (\ppt.2) \hskip5mm
\frac{{\CC_{(0)}}^{\ptr m}}{m!} \nts \frac{\Theta^h}{h!} \ = \ {g-m+h \choose h } \frac{{\CC_{(0)}}^{\ptr m-h}}{(m-h)!} \ ,
\quad m , \, h \, \in \, \Z \ . $
\eitem

We shall see that the tautological ring $ \, \RR(\CC) \, $ can be decomposed as the direct sum of \lq\lq remarkable" subspaces,
i.e. having a nice behavior under the Fourier transform and both products.
This decomposition, in particular, allows us to describe the structure of $ \, \RR(\CC) \, $ and to estimate its dimension.
To this purpose, we introduce a few notations and a technical lemma.

{\bf \noindent Notation \pzn. \hskip3mm}
We consider the family $ \, \JJ \, $ of non-ordered sets of strictly positive integers
$ \ I \, = \, \{i_1 , \, ... , \, i_d\} \ $ and we define

\bitem
$ d(I) \ = \ $ cardinality of $ \ I \ ; $

$ s(I) \ = \ i_1 + ... + i_d \ . $
\eitem

\noindent
The empty set $ \, \emptyset \, $ is also admitted.
We may note that \ $ I = \emptyset $ \ if and only if \ $ d(I) = 0 \, , $ \ if and only if \ $ s(I) = 0 \, . \ $
We shall also use the following notation: for a cycle $ \, \lambda_I \, $ we set
$$
\lambda_I^{[m]} \ := \ \lambda_I \ptr {\frac{\ssize{\CC_{(0)}}^{\ptr m}}{\dsize m!}} \, , \ m \in \Z \ ; \qquad
\Lambda_I \ := \ \left\langle \, \lambda_I^{[m]} \, \right\rangle_{m \in \N} \ ,
\tag{\pzn.1}
$$
i.e. $ \, \Lambda_I \, $ is the vector space (over $\Q ) $ generated by the $ \, \lambda_I^{[m]} \, . \ $
Given two integers \ $ n , \, s $ \ and cycles $ \, \lambda_I \, $ defined for all $ \, I \in \JJ \, $ satisfying
$ \, d(I) = n \, $ and $ \, s(I) = s \, , \ $ we set
$$
\Lambda_{[n;\,s]} \quad := \quad \sum_{I \, \vert \, d(I) = n , \ s(I) = s} \ \Lambda_I
$$
Furthermore, for the $ \, n $-dimensional component of the tautological ring we will use the notation
$ \ \RR_n \, = \, \RR(\CC) \cap \AA_n\big(X\big)_{\Q} \, . \ $

\bproc{Lemma \pzl}
Let $ \, \CC \, $ be a curve of genus $ \, g \, , \ $ let $ \, (X,\Theta) \, $ be its principally polarized Jacobian and
let $ \ \RR \, = \, \RR(\CC) \ $ denote its tautological ring.
Then there exist cycles \ $ \lambda_I \ \in \ \left[\AA_{d(I)}\big(X\big)_{\Q} \right]_{s(I)} \, , $
\ $ \ I \, = \, \{i_1 , \, ... , \, i_d\} \, \in \, \JJ \, , $ \ such that

\bitem
$ (\pzl.1) \qquad \dsize \lambda_I \nts \Theta \quad = \quad 0 \, ; $

$ (\pzl.2) \qquad \dsize
\lambda_I \quad = \quad \CC_{(i_1)} \ptr ... \ptr \CC_{(i_d)} \ , \qquad \modulo \quad
\sum_{n < d(I) } \quad \Lambda_{[n;\,s(I)]} $

$($in particular $ \ \lambda_{\emptyset} \, = \, \{o\} \, , \ $ the unit of \ $ \RR $ \ with respect to the Pontryagin product$);$

$ (\pzl.3) \qquad \dsize
\sum_{n \le t \, , \ s \in \N} \quad \Lambda_{[n;\,s]} \quad \supseteq \quad
{\dsize \sum_{n \le t}} \ \RR \cap \AA_n\big(X\big)_{\Q} \ , \qquad \forall \ t \in \N \ . $
\eitem
\eproc

Before giving the proof, we want to state a few straightforward considerations.
Consider a cycle $ \, \lambda_I \, $ as in lemma (\pzl) and let $ \, s = s(I) \, , \ d = d(I) \, ; \ $
note that $ \ \lambda_I^{[m]} \, \in \, \left[\AA_{d+m}\big(X\big)_{\Q} \right]_s \, . \ $
Then, assuming $ \ \lambda_I \ne 0 \, , \ $ (\ppz.1) to (\ppz.4) of proposition (\ppz) become

\bitem
(\pzm.1) \qquad $ s + 2 d \ \le \ g \ ; $

(\pzm.2) \qquad $ \dsize \lambda_I^{[m]} \, \ne \, 0 \qquad \Longleftrightarrow \qquad 0 \ \le \ m \ \le \ g\!-\!s\!-\!2d \, ; $

(\pzm.3) \qquad $ \dsize \FF\left(\lambda_I^{[m]}\right) \ = \ (-1)^m \, \lambda_I^{[g-s-2d-m]} \ ; $

(\pzm.4) \qquad $ \lambda_I^{[m]} \, \nts {\dsize \Theta^h \over h !} \ = \ {g-s-2d-m+h \choose h} \ \lambda_I^{[m-h]} \ , \quad
\forall \ m , \, h \, \in \, \Z \ . $
\eitem

\noindent
We may note that, in particular, the $ \, \Lambda_I $'s
are stable under the Fourier transform and under the intersection with the theta divisor,
moreover they are stable under the Pontryagin product with $ \, \CC_{(0)} \, $ by definition.
Furthermore, as the $ \, \Lambda_I $'s have a basis of pure cycles, they are also stable with respect to
pull-backs and push-forwards of multiplication maps by integers.
Thus, clearly the $ \, \Lambda_{[n;\,s]} $'s satisfy the same stability properties.

As a consequence of (\pzl.2), we are forced to define $ \ \lambda_{\{s\}} \, = \, \CC_{(s)} \, , \ $ where $ \, s \ge 1 \, $
(recall that the $ \, \lambda_I $'s are defined only for strictly positive $ \, I $'s).
Another consequence is that for $ \, a_1 , \, ..., \, a_n \, \in \, \N \, $ satisfying $ \ s \, := \, \sum a_i \, > \, g\!-\!2n \ $ we have
$$
\CC_{(a_1)} \ptr ... \ptr \CC_{(a_n)} \quad \in \quad \sum_{i=0}^{n-1} \ \Lambda_{[i,\, s]} \cap \left[\AA_n\big(X\big)_{\Q} \right]_s
\tag{\pzm.5}
$$
(we want to stress that the summation above is only up to $ \, n-1) \ $
indeed, assuming that all the $ \, a_i \, $ are nonzero (otherwise the claim is trivial),
one has $ \, \lambda_{\{a_1 , \, ..., \, a_n\}} = 0 \, $ by (\pzm.2), then (\pzm.5) follows from (\pzl.2).
In dimension 1, one finds the known fact that $ \ \CC_{(s)} \, = \, 0 \ $ for $ \ s \, > \, g - 2 \, . \ $
In dimension 2, we get $ \, \lambda_{\{i,\,j\}} = 0 \, $ for $ \, i+j \, \ge \, g-3 \, $
(assume $ i , \, j \ne 0 ), \, $ and, therefore, by (\pzl.2) and (\pzm.2),
$$
\eqalign{
& \CC_{(i)} \ptr \CC_{(j)} \ \in \ \Lambda_{[1,\, i+j]} \cap \left[\AA_2\big(X\big)_{\Q} \right]_{i+j}
\ = \ \Q \, \CC_{(i+j)} \ptr \CC_{(0)} \ , \quad i\!+\!j \, = \, g\!-\!3 \, ; \cr
& \CC_{(i)} \ptr \CC_{(j)} \ \in \ \Lambda_{[1,\, i+j]} \cap \left[\AA_2\big(X\big)_{\Q} \right]_{i+j}
\ = \ 0 \ , \quad i\!+\!j \, \ge \, g\!-\!2 \, ;
}
\tag{\pzm.6}
$$
More precisely, we will see that for $ \, i\!+\!j = g\!-\!3 \, , \ $ one has
$ \ \CC_{(i)} \ptr \CC_{(j)} \, = \, - {\ssize{g-1 \choose i+1}} \, \CC_{(g-3)} \ptr \CC_{(0)} $
\ (cf. \rij \ below).
In dimension 3, using that $ \, \lambda_{\{i,\,j,\,h\}} = 0 \, $ for $ \, i+j+h \, > \, g-6 \, $
(assume $ i , \, j , \, h $ strictly positive as in \pzn) and using (\pzm.5) it is easy to check that
$$
\eqalign{
& \CC_{(i)} \ptr \CC_{(j)} \ptr \CC_{(h)} \ \in \
  \left\langle \, \CC_{(r)} \ptr \CC_{(g-5-r)} \ptr \CC_{(0)} \, \right\rangle_{1 \le r \le g-5} \, , \
  i\!+\!j\!+\!h \, = \, g\!-\!5 \, ; \cr
& \CC_{(i)} \ptr \CC_{(j)} \ptr \CC_{(h)} \ \in \
  \Q \, \CC_{(g-4)} \ptr \CC_{(0)}^{\ptr 2} \, , \ i\!+\!j\!+\!h \, = \, g\!-\!4 \, ; \cr
& \CC_{(i)} \ptr \CC_{(j)} \ptr \CC_{(h)} \ = \ 0 \, , \ i\!+\!j\!+\!h \, \ge \, g\!-\!3 \ .
}
$$

\proof{\ (of lemma \pzl)}
We proceed by induction on dimension.
Fix $ \, \delta \ge 0 \, $ and assume the $ \, \lambda_I $'s are defined for all \ $ I \in \JJ $ \ satisfying
$ \ d(I) \le \delta \, . \ $
Let $ \ \Omega \, := \, \left\langle \, \lambda_I^{[m]} \, \right\rangle_{d(I) \le \delta , \, m \in \Z} \ $
and let
$$
Y \quad := \quad \CC_{i_1} \ptr ... \ptr \CC_{i_{\delta+1}} \ , \quad s \, := \, s(\{i_1, \, ..., \, i_{\delta+1}\}) \ .
$$
According with (\pzl.1) and (\pzl.2), we search, in particular, for a cycle
$ \, \lambda \, = \, \lambda_{\{i_1, \, ..., \, i_{\delta+1}\}} \, $
satisfying $ \ \lambda = Y \ \modulo \ \Omega \ $ and $ \ \lambda \nts \Theta = 0 \, . $
Finding such a cycle $ \, \lambda \, $ is equivalent to proving that there exists a cycle $ \ Z \, \in \, \Omega \ $ such that
$ \ Y \nts \Theta \, = \, Z \nts \Theta $ \ (define $ \, \lambda = Y\!-\!Z ). $
\ Recall that $ \ \Omega \ $ is stable under the Fourier transform and intersection with the theta divisor.

We distinguish two possibilities. First, by (\fut.3), assuming $ \ g-s-\delta-1 \, \le \, \delta \, , \ $
we have $ \ \dim \FF(Y) \, = \, g-s-\delta-1 \, \le \, \delta \, , \ $ so $ \ \FF(Y) \in \Omega \, . \ $
Thus $ \ Y \, = \, \pm \FF(\FF(Y)) \, \in \, \FF(\Omega) \, = \, \Omega \ $ and,
defining $ \ Z = Y \ $ (i.e. $ \, \lambda = 0), \, $ we are done.

We now consider the case when $ \ g-s-\delta-1 \, \ge \, \delta \, . \ $ We have
$$
\eqalign{
Y \nts \Theta \quad & = \quad (-1)^{g+s} \ \FF \big( \FF(Y \nts \Theta) \big) \cr
& = \quad (-1)^{g+s} \ \FF \left[ \sum_{j=0}^{\delta} (-1)^j \left(Y\nts \frac{\Theta^{j+1}}{j!} \right) \ptr
\frac{\Theta^{s+2\delta-j}}{(s+2\delta-j)!} \right] \cr
& = \quad \sum_{j=0}^{\delta} \ \left[ \FF\left(Y\nts \frac{\Theta^{j+1}}{j!} \right) \right] \nts
\frac{\Theta^{g-s-2\delta+j}}{(g\!-\!s\!-\!2\delta\!+\!j)!}
}
$$
where the first equality follows by (\fut.2),
for the second equality just apply (\okt.6) with $ \, W = Y \nts \Theta \, , \ $ and the third equality
follows by (\fut.1) and (\ppt.1).
Eventually, defining
$$
Z \quad := \quad {\dsize \sum_{j=0}^{\delta}} \
\left[ \FF\left(Y\nts \frac{\Theta^{j+1}}{j!} \right) \right] \nts \frac{\Theta^{g-s-2\delta+j-1}}{(g\!-\!s\!-\!2\delta\!+\!j)!}
$$
(note that we are under the hypothesis $ \ g\!-\!s\!-\!2\delta \ge 1 \, ), \ $
we have $ \ Y \nts \Theta \, = \, Z \nts \Theta \ $ as required.
We are left to prove that $ \, Z \in \Omega : \ $
from the induction hypothesis and (\pzl.3) we have
$ \, Y \nts \frac{\Theta^{j+1}}{j!} \, \in \, \RR_{\delta-1-j} \, \subseteq \, \Omega \, . \ $
As $ \ \Omega \ $ is stable under the Fourier transform and intersection with the theta divisor,
each term from the summation which defines $ \, Z \, $ is in $ \, \Omega \, , \ $ and we are done.

So far we have defined $ \, \lambda_I \, $ for $ \, d(I) \le \delta\!+\!1 \, $ and
such $ \, \lambda_I $'s satisfy (\pzl.1) and (\pzl.2).
We are left to prove (\pzl.3) with $ \, t = \delta\!+\!1 \, . \ $ This is clear. In fact,
$$
\sum_{n \le \delta+1 \, , \ s \in \N} \ \Lambda_{[n;\,s]} \quad \supseteq \quad
\big\langle \CC_{(a_1)} \ptr ... \ptr \CC_{(a_n)} \big\rangle_{a_i \ge 0 \, , \ n \le \delta+1} \quad = \quad
{\dsize \sum_{n \le \delta+1}} \ \RR_n
$$
where equality follows by the fact that $ \, \RR \, $ is generated by the classes $ \, \CC_{(i)} \, , \ $ see [Be3].
\qed

\bproc{Corollary \pzc}
The spaces $ \ \Lambda_{[n;\,s]} \ $ are independent.
\eproc

\proof{}
As the spaces $ \ \left[\AA_d\big(X\big)_{\Q} \right]_{s} \ $ are independent
and the $ \, \Lambda_{[n;\,s]} \, $ are generated by homogeneous elements, it suffices to prove
that does not exist a non-trivial relation
$ \ \sum \alpha_n \, = \, 0 \, , \ \alpha_n \in \Lambda_{[n;\,s_0]} \cap \left[\AA_{d_0}\big(X\big)_{\Q} \right]_{s_0} \, , \ $
where $ \, d_0 \, $ and $ \, s_0 \, $ are fixed (we may note that $ \, d_0 \ge n \, , \ $ unless $ \, \alpha_n \, $ is trivial).

By contradiction, assume there is such a non-trivial relation and observe that each $ \, \alpha_n \, $ must be of the form
$ \, \sum c_i \lambda_{I_i}^{[d_0-n]}  \, , \ $ where $ \, d(I_i) = n \, , \ c_i \in \Q \, . \ $
Let $ \, \rho \, $ be the operator defined by
$$
\rho_h (\alpha) \quad := \quad \left[\alpha \nts \Theta^h \right] \ptr {\CC_{(0)}}^{\ptr h} \ .
$$
By lemma (\pzl) and (\pzm.4), we have
$ \ \rho_h (\alpha_n) \, = \, 0 \ $ for $ \, h > d_0\!-\!n \, , $ \ and also
$ \ \rho_{d_0-n} (\alpha_n) \, = \, \frac{(d_0-n)!(g-s_0-2n)!}{(g-s_0-d_0-n)!} \, \alpha_n \, . $
\ Let now $ \ n' \, = \, \min \{\, n \, \vert \, \alpha_n \ne 0 \} \, . \ $
Then,
$$
0 \quad = \quad \rho_{d_0-n'}\big(\sum \alpha_n \big) \quad = \quad \frac{(d_0-n)!(g-s_0-2n')!}{(g-s_0-d_0-n')!} \, \alpha_{n'} \quad \ne \quad 0 \ .
$$
\ This is a contradiction.
\qed

\bproc{Corollary \pzr}
Let $ \, \delta \, $ and $ \, s \, $ be integers.
The restriction to $ \ \left[\AA_\delta(X)_\Q\right]_s \cap {\dsize\sum_{n<\delta}} \Lambda_{[n;\,s]} \ $
of the intersection with the Theta divisor map is injective.
\eproc

\proof{}
By corollary (\pzc) and the fact that each $ \, \Lambda_{[n;\,s]} \, $ is stable under the intersection with the Theta divisor,
it suffices to prove that the intersection with the Theta divisor map is injective on
$ \ \left[\AA_\delta(X)_\Q\right]_s \cap \Lambda_{[\nu;\,s]} \ $ (for a fixed $ \, \nu < \delta ). \ $
Consider a non-zero cycle $ \ \alpha \, \in \, \left[\AA_\delta(X)_\Q\right]_s \cap \Lambda_{[\nu;\,s]} \, . \ $
Thus, $ \ \alpha \, = \, \sum c_i \lambda_{I_i}^{[\delta-\nu]} \, , \ $
where $ \, d(I_i) = \nu \, $ and $ \, s(I_i) = s \, . $
\ As $ \, \alpha \ne 0 \, , \ $ from (\pzm.2), we have $ \ \delta - \nu \, \le \, g-s-2\nu \, . \ $
On the other hand, from (\pzm.4) and (\pzn), we have
$ \ (\alpha \nts \Theta)\ptr \Gamma \, = \, (g-s-2\nu-\delta+\nu+1)(\delta-\nu) \pt \alpha \, . $
\ As the coefficient is non-zero, also $ \, \alpha \nts \Theta \, $ cannot be trivial, and we are done.
\qed

\bproc{Corollary \pzd}
The set of $ \ \lambda $'s from lemma $(\pzl)$ is unique.
\eproc

\proof{}
Consider two such sets $ \, \{\lambda_I\} \, $ and $ \, \{\lambda'_I\} \, $ and let
$ \, J \, $ satisfy $ \, \lambda_J \ne \lambda'_J \, $ with $ \, \delta := d(J) \, $ minimal
(so that $ \, \Lambda_{[n;\,s]} = \Lambda'_{[n;\,s]} \, , \ \forall \ n < \delta). $
\ Then, from (\pzl.1) and (\pzl.2), we have $ \ (\lambda_J - \lambda'_J) \nts \Theta \, = \, 0 \ $ and
$ \ \lambda_J - \lambda'_J \, \in \, {\dsize\sum_{n<\delta}} \Lambda_{[n;\,s(J)]} \, . \ $
This contradicts corollary (\pzr).
\qed

\bproc{Corollary \pzv}
If $ \ \CC_{i_1} \ptr ... \ptr \CC_{i_d} \, = 0 \, , \ $ then $ \ \lambda_{\{i_1,\, ...,\, i_d\}} \, = \, 0 \ . $
\eproc

\proof{}
This is clear by the previous corollary, by the proof of lemma (\pzl), and by the fact that the new set of
$ \, \lambda $'s, obtained substituting
$ \, \lambda_{\{i_1,\, ...,\, i_d\}} \, $ with $ \, 0 \, , \ $ does satisfy (\pzl.1), (\pzl.2) and (\pzl.3).
\qed

\

Lemma (\pzl) and its corollaries  are a tool to give an explicit description of the tautological ring $ \, \RR(\CC) \, . \ $
In fact, using the $ \, \lambda_I^{[m]}  \, $ we get a nice presentation of the tautological ring,
which we will write down explicitly for $ \, g \le 8 \, , \ $ also listing all the possibilities that may occur and giving a model for the
bi-graded two-products ring $ \, \RR(\CC) \, . \ $
To this purpose we need to discuss the matrices of the algebraic structure of our ring $ \, \RR(\CC) \, ; \ $
namely the matrices representing $ \, \FF \, , \ $ the Pontryagin product
and the intersection product, corresponding to our basis of $ \, \lambda_I^{[m]} $ (with non-trivial $ \, \lambda_I $'s).

The reader is advised to look at pictures 1 to 6 below, at section 6 (clearly, one can draw similar pictures for any $ \, g \, ). \ $

\rmk{Remark \rbp}
First, according with (\pzm.2), non triviality of $ \, \lambda_I^{[m]} \, $ depends only on $ \, I \, $
(consider $ m $ in the range from $ 0 $ to $ g-s(I)-2d(I)) . \ $
As for the matrix representing $ \, \FF \, , \ $ it is clear from (\pzm.3) that it is a block matrix
$ \ \mmm{B_1}{\ddots}{B_k} \, , \ $ where $ \, B_i \, $ is the matrix with $ \, \{ 1 , \, -1 , \, ... \} \, $ on
the secondary diagonal (starting with +1 at the left-bottom corner) and zero elsewhere,
and where each block corresponds to the elements $ \, \lambda_I^{[0]} , \, ..., \, \lambda_I^{[g-s(I)-2d(I)]} . \ $
In order to write down explicitly the matrices representing the Pontryagin product and the intersection product,
it is required a further computation, at least for higher genus cases: the definition itself of the $ \ \lambda_{I}^{[m]} $'s
takes care of the Pontryagin product of an element from the first column with another element from the picture
while (\pzm.4) takes care of intersection products of elements from the first column with other elements from the picture.
To be very explicit, this amounts to the following
$$
\eqalign{
& \lambda_{\emptyset}^{[k]} \, \ptr \, \lambda_I^{[m]} \quad = \quad {\ssize{m+k \choose m}} \ \lambda_I^{[m+k]}\cr
& \lambda_{\emptyset}^{[k]} \ \nts \ \lambda_I^{[m]} \quad = \quad {\ssize{2g-m-k-s-2d \choose g-k}} \ \lambda_I^{[m-g+k]}
}
\tag{\mpt}
$$
It is also clear that we get zero whenever a product land outside the picture
(this event can be easily checked thanks to the fact that the Pontryagin product is dimensionally homogeneous,
the intersection product is co-dimensionally homogeneous and both products are homogeneous with respect to Beauville's graduation).
Furthermore, in sight of the straightforward identities
$$
\eqalign{
& \lambda_I^{[m]} \, \ptr \, \lambda_J^{[h]} \quad = \quad {\ssize{m+h \choose m}} \
  \lambda_{\emptyset}^{[m+h]} \, \ptr \, \lambda_I \, \ptr \, \lambda_J \cr
& \lambda_I^{[t]} \ \nts \ \lambda_J^{[k]} \quad = \quad {\ssize{i'-t+j'-k \choose i'-t}} \
  \lambda_{\emptyset}^{[g-i'+t-j'+k]} \, \nts \, \lambda_I^{[i']} \, \nts \, \lambda_J^{[j']}
}
\tag{\mpg}
$$
where $ \, i' := g-s(I)-2d(I) \, $ and $ \, j' := g-s(J)-2d(J) \, , \ $ the products of the form
$ \ \lambda_I \ptr \lambda_J $ \ and $ \ \lambda_I^{[i']} \nts \lambda_J^{[j']} \ $
determine all other products.
We shall refer to these products as basic products
(to avoid long lists we shall limit ourselves to write down only basic products when describing intersection and Pontryagin matrices).
We may also note that the two equalities above are $ \FF $-dual to each other (under $ \, t = i'-m \, , \ k = j'-h \pt ), \ $
in particular basic intersection products are the $ \FF $-dual of basic Pontryagin products:
$$
\lambda_I^{[i']} \, \nts \, \lambda_J^{[j']} \quad = \quad \FF \pt \left( \lambda_I \ptr \lambda_J \right)
\tag{\mpb}
$$

We end this section with further technical results. First, our goal is to compute explicitly the $ \ \lambda_{\{i,\, j\}} $'s.
From lemma (\pzl), in particular from (\pzl.2), we have $ \ \lambda_{\{i,\, j\}} = 0 \ $ for $ \ i+j \, \ge \, g-3 \, . $
\ Furthermore, by (\pzm.6), $ \ \lambda_{\{i,\, j\}} \ $ must be equal to
$ \ \CC_{(i)} \ptr \CC_{(j)} + \alpha \, \CC_{(0)} \ptr \CC_{(i+j)} $ \ for some $ \, \alpha \in \Q \, , $ \ and
in order to satisfy (\pzl.1), \, i.e. $ \lambda_{\{i,\, j\}} \nts \Theta = 0 \, , \ $
the coefficient $ \, \alpha \, $ must satisfy
$ \ \left(\CC_{(i)} \ptr \CC_{(j)}\right) \nts \Theta = - \alpha ({g\!-\!i\!-\!j\!-\!2}) \CC_{(i+j)} $
\ (see also \pzm.4).
Lemma (\cij) below computes the intersection at the left hand side and eventually determines the
coefficient $ \, \alpha \, . \ $

\bproc{Lemma \cij}
Let $ \, \CC \, $ be a curve of genus $ \, g \, , \ $ and let $ \ i , \, j \, \ge \, 1 \, . \ $ Then
$$
\left( \CC_{(i)} \ptr \CC_{(j)} \right) \nts \Theta \quad = \quad - {i+j+2 \choose i+1} \, \CC_{(i+j)} \ ,
$$
\eproc

As a consequence, we get $ \ \alpha \, = \frac{{i+j+2\choose i+1}}{g\!-\!i\!-\!j\!-\!2} \ $ (for $ \, i+j \le g-3 ). \ $ Thus
$$
\lambda_{\{i,\, j\}} \quad = \quad \left\{
\ \eqalign{
& 0 \ , \hskip60mm i+j \, \ge \, g-3 \ ; \cr
& \CC_{(i)} \ptr \CC_{(j)} + \frac{{i+j+2\choose i+1}}{g\!-\!i\!-\!j\!-\!2} \, \CC_{(i+j)} \ptr \CC_{(0)} \ , \qquad i+j \, \le \, g-3 \ .
}
\right.
\tag{\rij}
$$
In particular, for $ i+j = g-3, $ we have $ \CC_{(i)} \ptr \CC_{(j)} = - {g-1 \choose i+1} \CC_{(g-3)} \ptr \CC_{(0)} . \ $
It should be observed that $ (\rij) $ can be read as a result that regards basic products:
$$
\lambda_{\{i\}} \ptr \lambda_{\{j\}} \quad = \quad
\lambda_{\{i,\, j\}} \ - \ \frac{{i+j+2\choose i+1}}{g\!-\!i\!-\!j\!-\!2} \, \lambda_{\{i+j\}}^{[1]} \ , \qquad i+j \, \le \, g-3
\tag{\rijb}
$$
(and such basic product is trivial for $ \, i+j \ge g-2 ). $

\bproc{Corollary \cvn}
If $ \ C_{(i)} = 0 \, , \ $ then $ \ C_{(s)} = 0 \, $ for all $ \, s \ge i \, . $
\eproc

\proof{}
By the previous lemma, the cycle $ \ C_{(s)} \ $ is a nontrivial multiple of
$ \ \left( \CC_{(i)} \ptr \CC_{(s-i)} \right) \nts \Theta \, , $ \ for $ \ s > i \, . $
\qed

\proof{ (of lemma \cij)}
Consider the composition
$$
u : \ \CC \times \CC \ \blr{\lra}{\phi} \ J(\CC) \times J(\CC) \ \blr{\lra}{(\cdot m ; \cdot n)} \ J(\CC) \times J(\CC) \ \blr{\lra}{\sigma} \ J(\CC)
$$
where $ \, \phi \, $ is the natural inclusion, $ \ (\cdot m ; \cdot n) \ $ denotes the multiplication map
(by $ \, m \, $ on the first factor and by $ \, n \, $ on the second factor) and $ \ \sigma \, $ denotes the sum map.
As
$ \, \sigma^{\star} \Theta \, = \, \Theta \times J(\CC) + J(\CC) \times \Theta - \PP \, , \ $ we have
$ \, (\cdot m ; \cdot n)^{\star} \sigma^{\star} \Theta \, = \, m^2 \Theta \times J(\CC) + n^2 J(\CC) \times \Theta - mn\PP \, , $
\ therefore $ \ u^{\star}\Theta \, = \, (m^2g+mn)\{o\} \times \CC + (n^2g+mn) \CC \times \{o\} - mn \Delta \, , \ $
where $ \, \Delta \, $ denotes the diagonal of $ \, \CC \times \CC \, . \ $ Taking $ \, u_{\star} \, $ we obtain
$$
\eqalign{
& \big(\mls{m}\CC \ptr \mls{n}\CC \big) \nts \Theta \quad \cr
& = \quad (m^2g+mn)\mls{n}\CC + (n^2g+mn)\mls{m}\CC - mn \mls{m+n}\CC
}
\tag{\cij.1}
$$
Let now $ \ \sqb{s}{n} \, := \, \dsize \frac{(-1)^{s+n} {s+2 \choose n}}{(s+2)!} \, . \ $
Because of the identity (see [Ma], proposition 9)
$$
\CC_{(s)} \quad = \quad \sum_{n=1}^{s+2} \sqb{s}{n} \mls{n}\CC \ , \qquad
\trm{modulo } \ \sum_{k>s} \, \big\langle \CC_{(k)} \big\rangle_\Q \ ,
$$
we have
$$
\eqalign{
& \left( \CC_{(i)} \ptr \CC_{(j)} \right) \nts \Theta \quad \cr
& = \quad
\left( \left[\sum_{m=1}^{i+2} \sqb{i}{m} \mls{m}\CC \right] \ptr
       \left[\sum_{n=1}^{j+2} \sqb{j}{n} \mls{n}\CC \right] \right) \nts \Theta \cr
& = \quad
\sum_{m=1}^{i+2} \sum_{n=1}^{j+2} \sqb{i}{m} \sqb{j}{n} \big( \mls{m}\CC \ptr \mls{n}\CC \big) \nts \Theta
}
\tag{\cij.2}
$$
modulo $ \ \sum_{k>i+j} \, \big\langle \CC_{(k)} \big\rangle_\Q \, . \ $
Now, plugging (\cij.1) into the right hand side of (\cij.2), expanding $ \, \CC \, $ as $ \, \sum \CC_{(s)} \, , \ $
performing the substitutions $ \ \mls{n} \CC_{(s)} = n^{s+2} \CC_{(s)} \ $
and omitting terms in $ \ \sum_{k>i+j} \, \big\langle \CC_{(k)} \big\rangle_\Q \ $
we obtain the expression
$$
\sum_{t=0}^{i+j} \ \xi_{i,\,j}^{[t]} \ \CC_{(t)} \ ,
$$
where
$$
\xi_{i,\,j}^{[t]} \ := \ {\dsize \sum_{m=1}^{i+2} \sum_{n=1}^{j+2}}
\sqb{i}{m} \cdot \sqb{j}{n} \cdot
\left[\big((m^2g+mn)n^{t+2}+(n^2g+mn)m^{t+2}-mn(m+n)^{t+2}\big) \right]
$$
As a consequence
$$
\big( \CC_{(i)} \ptr \CC_{(j)} \big) \nts \Theta \quad = \quad \xi_{i,\,j}^{[i+j]} \ \CC_{(i+j)}
\tag{\cij.3}
$$
(as clearly expected, we have $ \, \xi_{i,\,j}^{[t]} \, = \, 0 \, $ for $ \, t < i\!+\!j \, ). \ $
Then, by a straightforward computation,
$$
\xi_{i,\,j}^{[i+j]} \quad = \quad \left\{
\eqalign{
& (2g-2) \ , \hskip20mm i=j=0 \cr
& (g-j-2) \ , \hskip15.7mm i=0, \ j\ge 1 \cr
& - {i+j+2 \choose i+1} \ , \hskip10mm i\ge 1, \ j\ge 1
}
\right.
\tag{\cij.4}
$$
and we are done.
\qed

It is worth to note that (\cij.3) and (\cij.4) also give
$ \ \CC_{(0)}^{\ptr 2} \nts \Theta \, = \, (2g\!-\!2) \CC_{(0)} \ $ and
$ \ \big( \CC_{(s)} \ptr \CC_{(0)} \big) \nts \Theta \, = \, (g\!-\!s\!-\!2) \CC_{(s)} \, $ for $ \, s \ge 1 \, , \ $
as already known from (\tgf.2) and (\pzm.4) respectively.

\bproc{Lemma \chij}
Let $ \, \CC \, $ be a curve of genus $ \, g \, , \ $ and let $ \ h , \, i , \, j \, \ge \, 1 \, . \ $ Then
$$
\eqalign{
\left( \CC_{(h)} \ptr \CC_{(i)} \ptr \CC_{(j)} \right) \nts \Theta \qquad = \qquad
& - {\ssize {i+j+2 \choose i+1}} \, \CC_{(h)} \ptr \CC_{(i+j)} \cr
& - {\ssize {j+h+2 \choose j+1}} \, \CC_{(i)} \ptr \CC_{(j+h)} \cr
& - {\ssize {h+i+2 \choose h+1}} \, \CC_{(j)} \ptr \CC_{(h+i)}
}
$$
\eproc

\proof{}
Considering the composition
$$
u : \ \CC \times \CC \times \CC \ \blr{\lra}{\vec\phi} \ J(\CC) \times J(\CC) \times J(\CC) \
\blr{\lra}{(\cdot l ; \cdot m ; \cdot n)} \ J(\CC) \times J(\CC) \times J(\CC) \ \blr{\lra}{\sigma} \ J(\CC)
$$
and proceeding as in the proof of lemma (\cij), we obtain
$$
\eqalign{
& \big(\mls{l}\CC \ptr \mls{m}\CC \ptr \mls{n}\CC \big) \nts \Theta \quad \cr
= \quad \cr & \eqalign{   & (n^2g+nl+nm) \; \mls{l}\CC \ptr \mls{m}\CC \cr
                        + \ & (l^2g+lm+ln) \; \mls{m}\CC \ptr \mls{n}\CC \cr
                        + \ & (m^2g+mn+ml) \; \mls{n}\CC \ptr \mls{l}\CC \cr
                        - \ & (lm) \; \mls{l+m}\CC \ptr \mls{n}\CC \cr
                        - \ & (mn) \; \mls{m+n}\CC \ptr \mls{l}\CC \cr
                        - \ & (nl) \; \mls{n+l}\CC \ptr \mls{m}\CC
}
}
\tag{\chij.1}
$$
On the other hand we also have the analogous of (\cij.2),
$$
\eqalign{
& \left( \CC_{(h)} \ptr \CC_{(i)} \ptr \CC_{(j)} \right) \nts \Theta \quad \cr \cr
& = \quad
\sum_{l=1}^{h+2} \sum_{m=1}^{i+2} \sum_{n=1}^{j+2} \sqb{h}{l} \sqb{i}{m} \sqb{j}{n}
\big( \mls{l}\CC \ptr \mls{m}\CC \ptr \mls{n}\CC \big) \nts \Theta
}
\tag{\chij.2}
$$
modulo $ \ \big[A_2(J(\CC))_\Q\big]_{>\, h+i+j} \, . \ $
Then, proceeding as in the previous lemma, and precisely plugging (\chij.1) into the right hand side of (\chij.2),
expanding $ \, \CC \, $ as $ \, \sum \CC_{(s)} \, , \ $
performing the substitutions $ \ \mls{n} \CC_{(s)} = n^{s+2} \CC_{(s)} \ $
and omitting terms in $ \ \big[A_2(J(\CC))_\Q\big]_{>\, h+i+j} \, , $ \ we obtain an expression of the following type:
$$
\sum_{r+t\,\le \,h+i+j} \ \xi_{h,\,i,\,j}^{[r,\,t]} \ \CC_{(r)} \ptr \CC_{(t)} \ .
$$
It is a very long still straightforward the computation of the coefficients
$ \, \xi_{h,\,i,\,j}^{[r,\,t]}  \, $ to check that they are the expected ones.
\qed

\newpage

\prg{5}{Algebraic models for $ \ \RR(\CC) \, $}

So far we have given generators for $ \, \RR(\CC) \, $ and we have collected enough pieces of information about their behavior
(under products, the Fourier transform, $ \mus{n} $ and $ \mls{n} $).
A priori it is not clear whether the jumble of formulas at our disposal hydes some further
relation (indeed, for example, we shall see that, for $ g=7 $, we must have $ \, \lambda_{\{2,\,1\}} = 0 ). $
\ Answer to this question means to determine admissible algebraic models.
This reduces to quite elementary algebra, the results are collected in this section.

\rmk{Definition \ami}
Let $ \, g \ge 1 \, $ and let $ \, \JJ \, $ be the family from notation (\pzn).
A set $ \, \AA \, $ of multi indexes is said {\it g-admissible} if
$$
\eqalign{
& \emptyset \in \AA \subseteq \JJ \ ; \qquad s(I) + 2d(I) \, \le \, g \, , \ \forall \ I \in \AA \ ; \cr
& I \in \AA \quad \Rightarrow \quad J \in \AA \ , \quad \forall \ J \, \subseteq \, I
}
$$
We consider symbols $ \, \lambda_I^{[m]} \, . \ $ We also consider symbols $ \mls{n} $ and $ \mus{n} $
operating on the $ \, \lambda_I^{[m]} \, $ as the multiplication by
$ \, n^{2d(I)+2m+s(I)} \, $ and $ \, n^{2g-2d(I)-2m-s(I)} \, $ respectively.

\bproc{Proposition \als}
Consider an integer $ \, g \ge 1 \, , \ $ a set $ \ \AA \ $ of g-admissible multi indexes and the vector space
$$
\trm{\bf V} \quad := \quad \bigoplus_{I \in A} \ \bigoplus_{m \, = \, 0}^{g-s(I)-2d(I)} \ \Q \ \lambda_I^{[m]} \ .
$$
endowed with the two graduations \lq\lq$ \pt s(I) \pt $" and \lq\lq$ \pt d(I) + m \pt $".
\ For $ \, I \in \JJ \, $ and $ \, m \in \Z \, $ also define $ \, \lambda_I^{[m]} = 0 \, $ if either $ \, I \not\in \AA \, $
or $ \, m \, $ is not in the range $ \ [ \ 0 \ , \ g-s(I)-2d(I) \, ] \, . $ \ Let $ \, \lambda_I := \lambda_I^{[0]} \, , $
\ consider stable $($with respect to both graduations$)$ \, {\bf V}-valued commutative \lq\lq basic" products
$ \ \lambda_I \ptr \lambda_J $ \ satisfying $ \ \lambda_{\emptyset} \ptr \lambda_I \equiv \lambda_I \ $ and define

$ (\als.1) \qquad
\lambda_{\emptyset}^{[k]} \, \ptr \, \lambda_I^{[m]} \quad = \quad {m+k \choose m} \ \lambda_I^{[m+k]} $

$ (\als.2) \qquad
\lambda_I^{[m]} \ \ptr \ \lambda_J^{[h]} \quad = \quad {m+h \choose m} \
\lambda_{\emptyset}^{[m+h]} \ \ptr \ \big(\lambda_I \ \ptr \ \lambda_J \big) $

Assuming $($associativity for basic products$)$
$$
(\lambda_I \ptr \lambda_J) \ptr \lambda_R \, \equiv \, \lambda_I \ptr (\lambda_J \ptr \lambda_R) \ ,
\tag{\als.3}
$$
there exist products \lq\lq$\ptr$" and \lq\lq$\nts$" and a endomorphism $ \, \FF \, $ satisfying
$$
(\fut.1), \ (\fut.2), \ (\fut.3), \ (\fut.4), \ (\pzm.1), \ (\pzm.2), \ (\pzm.3), \ (\mpt), \ (\mpg).
$$
Furthermore, such \lq\lq$\ptr$", \lq\lq$\nts$" and $ \FF $ are unique.
\eproc

\newpage

\proof{}
First, it is easy to check that (\als.1) and (\als.2) are compatible and that {\bf V} is a ring under \lq\lq$\ptr$".
Then, clearly one has to define the endomorphism $ \, \FF \, $ via (\pzm.2), \ i.e. \
$ \, \FF\left(\lambda_{I}^{[m]}\right) \ = \ (-1)^m \, \lambda_{I}^{[g-s(I)-2d(I)-m]} \, $
and \ \lq\lq $\nts$" via the $ \, \FF $-dual of \lq\lq$\ptr$", \ i.e. \ $ \lambda_{I}^{[m]} \nts \lambda_{J}^{[h]} \, = $
$ \, (-1)^{s(I)+s(J)} \FF \big[ \FF(\lambda_{I}^{[m]}) \ptr \FF(\lambda_{J}^{[h]}) \big] \, . \ $
Finally, is again straightforward to check that {\bf V} is a ring under \lq\lq$\nts$" and that all mentioned properties hold.
Uniqueness of such \lq\lq$\ptr$", \lq\lq$\nts$" and $ \FF $ is also clear.
\qed

This proposition leads to two heuristic considerations about models for $ \, \RR(\CC) \, : \ $ for the first one,
once basic Pontryagin products are given the structure of the model is uniquely determined; the second consideration is
that checking the existence of a model reduces to checking associativity for basic products.
The first non-trivial example of such verification is $ \, \big(\lambda_{\{1\}} \ptr \lambda_{\{1\}}\big) \ptr \lambda_{\{2\}} \,
= \, \lambda_{\{1\}} \ptr \big( \lambda_{\{1\}} \ptr \lambda_{\{2\}} \big) \, $ and it occurs only for $ \, g \ge 8 \, . $

Further geometrical relations, such as those one given by lemma (\cij) and lemma (\chij), give rise to further restrictions for
admissible models.

\rmk{Remark \alsb}
In terms of models, lemma (\cij) reduces to a condition on basic Pontryagin products:
namely it reduces to the identity (\rijb).
Note that by the construction of models, namely for reason of graduation and reason of $ g $-admissibility of multi indexes,
for $ \, i+j \, \ge \, g-2 \, $ one must have triviality for the basic products $ \, \lambda_{\{i\}} \ptr \lambda_{\{j\}} \, . \ $
We may also note that corollary (\cvn) can be stated as
$$
\{i\} \ \in \ \AA \qquad \Rightarrow \qquad \{j\} \ \in \ \AA \ , \qquad \forall \ j \, \le \, i
$$
(by the very definition of the model {\bf V}, clearly $ \ \{i\} \in \AA \, $ if and only if
$ \, \lambda_{\{i\}} \ne 0 \, ). $

Meanwhile lemma (\cij) reduces essentially to a condition on basic products, lemma (\chij) reduces to a more complicated requirement
and give rise to restrictions about admissible models.
In order to give a very explicit description for the possibilities that may occur for $ \, \RR(\CC) \, $
for $ \, g \le 8 \, , \ $ we shall use the following.

\rmk{Remark \alsc}
Assume $ \, g \ge 6 \, . \ $ By lemma (\chij) and formulae (\rijb), we must have
$$
\eqalign{
\left( \CC_{(1)} \ptr \CC_{(1)} \ptr \CC_{(1)} \right) \nts \Theta \qquad
& = \qquad - 18 \, \CC_{(1)} \ptr \CC_{(2)} \cr
& = \qquad - 18 \, \lambda_{\{2,\,1\}} \ + \ {180 \over g-5} \, \lambda_{\{3\}}^{[1]}
}
$$
and, for $ \, g \ge 7 \, , \ $ we have
$$
\eqalign{
\left( \CC_{(1)} \ptr \CC_{(1)} \ptr \CC_{(2)} \right) \nts \Theta \qquad
& = \qquad - 6 \, \CC_{(2)} \ptr \CC_{(2)} \ - \ 20 \, \CC_{(1)} \ptr \CC_{(3)} \cr
& = \qquad - 6 \, \lambda_{\{2,\,2\}} \ - 20 \, \lambda_{\{1,\,3\}} \ + \ {420 \over g-6} \, \lambda_{\{4\}}^{[1]}
}
$$

\newpage

\prg{6}{$ \RR(\CC) \, $ in low genus and other examples}

In this section we describe the admissible models for $ \, \RR(\CC) \, $ for $ \, g \le 8 \, . \ $
To this purpose, as generators of $ \, \RR(\CC) \, , \ $ we consider the $ \, \lambda_{I}^{[m]} \, . \ $
For $ \, g \le 5 \, $ the situation is very simple because only the
$ \, \lambda_{\emptyset}^{[m]} \, $ (that are the powers of the theta divisor)
and the $ \, \lambda_{\{s\}}^{[m]} \, $ may occur.
In fact, $ \, \lambda_{\{i,\,j\}} \, $ may occur only for $ \, g \ge i+j+2d \ge 6 \, . \ $
For the sake of completeness we draw the picture also for this cases
(anyway, for $ \, g \le 2 \, , \ \RR(\CC) $ \, has dimension $ \, g+1 \, $ and it is just generated by the powers of the theta divisor).
Once for all we recall that $ \ \Gamma = {\Theta^{g-1}\over (g-1)!} = \CC_{(0)} \ $ and
$ \ \lambda_{\{s\}} \, = \, \CC_{(s)} \, , \ s \ge 1 \, . \ $
We also observe the following.

\rmk{Remark \rbe}
The matrix representing the Fourier transform $ \, \FF \, $ is as described in remark (\rbp).
Therefore, in the pictures below, $ \, \FF \, $ acts as the vertical symmetry (up to a sign, see either (\pzm.3) or remark (\rbp)).
The intersection with the Theta divisor and the Pontryagin product with $ \, \CC_{(0)} \, $ correspond respectively
to a step downwards and a step upwards (up to a non-zero coefficient, see (\pzm.4) and the
definition of $ \, \lambda_{I}^{[m]} \, $ from notation (\pzn)).
The columns are independent (see corollary (\pzc)), and non-trivial columns have non-trivial elements by (\pzm.2).
Thus, the elements from non-trivial columns are a basis of $ \, \RR(\CC) \, . \ $
As for the matrices representing the Pontryagin and the intersection products (again, see remark (\rbp)),
they are governed by basic products, (\mpt), (\mpg) and (\mpb).

\vskip4mm

Case $ \, g = 3 \, . \ $

Other than the powers of the theta divisor (namely the column of $ \, \lambda_{\emptyset} ), \, $
by property (\csp) only $ \, \lambda_{\{1\}} \ \left(=\CC_{(1)}\right) \, $ may occur,
and it is nontrivial for the general Jacobian by a result of Ceresa (see [Ce]).
The picture is very simple, it is

\vskip-3mm$ \boxed{ \CD   
\hskip8mm J(\CC) \hskip8mm @. \hskip30mm \\
\Theta \\
\Gamma                     @. \lambda_{\{1\}}  \\
\{ \, o \, \}
\endCD } $ \vskip-3mm     

\vskip-10mm \rightline{picture (1)} \vskip6mm

\noindent
and, in particular, the admissible dimensions for $ \ \RR(\CC) \ $ are 4 and 5,
corresponding respectively to cases where $ \, \CC_{(1)} = 0 \, $ and $ \, \CC_{(1)} \ne 0 \, . $

\vskip4mm

Case $ \, g = 4 \, . \ $

Again by property (\csp), $ \, \CC_{(s)} \, = 0 \, $ for $ \, s \ge 2 \, $
and the situation is similar to the previous one, the picture is

\vskip-3mm$ \boxed{ \CD   
\hskip8mm J(\CC) \hskip8mm @. \hskip30mm          \\
\Theta                                            \\
\frac{\Gamma^{\ptr 2}}{2!} @. \lambda_{\{1\}}^{[1]} \\
\Gamma                     @. \lambda_{\{1\}}     \\
\{\,o\,\}
\endCD } $ \vskip-3mm     

\vskip-10mm \rightline{picture (2)} \vskip6mm

\noindent
so, the admissible dimensions for $ \, \RR(\CC) \, $ are 5 and 7,
corresponding respectively to cases where $ \, \CC_{(1)} = 0 \, $ and $ \, \CC_{(1)} \ne 0 \, $
(recall that nontrivial columns have nontrivial elements).

\vskip4mm

Case $ \, g = 5 \, . \ $

Again, by (\csp), we have $ \, \CC_{(s)} \, = 0 \, $ for $ \, s \ge 3 \, , \ $ thus the picture is as follows

\vskip-3mm$ \boxed{ \CD   
\hskip8mm J(\CC) \hskip8mm @. \hskip30mm          @. \hskip30mm              \\
\Theta                                                                       \\
\frac{\Gamma^{\ptr 3}}{3!} @. \lambda_{\{1\}}^{[2]}                          \\
\frac{\Gamma^{\ptr 2}}{2!} @. \lambda_{\{1\}}^{[1]} @. \lambda_{\{2\}}^{[1]} \\
\Gamma                     @. \lambda_{\{1\}}     @. \lambda_{\{2\}}         \\
\{\,o\,\}
\endCD } $ \vskip-3mm     

\vskip-10mm \rightline{picture (3)} \vskip6mm

\noindent
Then, the admissible dimensions for $ \, \RR(\CC) \, $ are 6, 9 and 11.
In case $ \, \CC_{(1)} = 0 \, , \ $ we must have also $ \, \CC_{(2)} = 0 \, $ by corollary (\cvn),
and the tautological ring has dimension 6.
Furthermore, assuming $ \, \CC_{(1)} \ne 0 \, $ and $ \, \CC_{(2)} = 0 \, , \ $ the tautological ring $ \, \RR(\CC) \, $ has dimension 9.
In case $ \, \CC_{(1)} \ne 0 \, $ and $ \, \CC_{(2)} \ne 0 \, $
(I don't know whether this case can actually occur), $ \, \RR(\CC) \, $ has dimension 11.
So, the possibilities that may occur are the following:

\bitem
(a) \quad $ \CC_{(1)} = 0                          \qquad (\Rightarrow \ \CC_{(2)} = 0 \, , \ \dim \RR(\CC) = 6) $

(b) \quad $ \CC_{(1)} \ne 0                        \qquad (\Rightarrow \ \CC_{(2)} = 0 \, , \ \dim \RR(\CC) = 9) $

(c) \quad $ \CC_{(1)} \ne 0 \, , \ \CC_{(2)} \ne 0 \qquad (\Rightarrow \ \dim \RR(\CC) = 11) $
\eitem

In view of remark (\rbe), as to describe the structure of $ \, \RR(\CC) \, , \ $
we are left to write down basic Pontryagin products.
Reminding that $ \, \lambda_{\emptyset} \ptr \lambda_I \equiv \lambda_I \, , \ $
we are only left to compute $ \, {\lambda_{\{1\}}}^{\ptr 2} \, $ in case (c): \ in view of (\rijb) we have
$$
{\lambda_{\{1\}}}^{\ptr 2} \ = \ -6 \, \lambda_{\{2\}}^{[1]} \ , \qquad \trm{case (c)}.
$$
Explicitly, in view of (\mpb), the corresponding basic intersection product is
$ \, \left(\lambda_{\{1\}}^{[2]}\right)^2 \ = \ 6 \, \lambda_{\{2\}} \, . $

\newpage

Case $ \, g = 6 \, . \ $

Again by property (\csp), we have $ \, \CC_{(s)} \, = 0 \, $ for $ \, s \ge 3 \, , \ $
the admissible dimensions are 7, 11, 12, 14, 15 and the picture is the following

\vskip-3mm$ \boxed{ \CD   
\hskip8mm J(\CC) \hskip8mm @. \hskip25mm          @. \hskip25mm          @. \hskip25mm           \\
\Theta                                                                                           \\
\frac{\Gamma^{\ptr 4}}{4!} @. \lambda_{\{1\}}^{[3]}                                                \\
\frac{\Gamma^{\ptr 3}}{3!} @. \lambda_{\{1\}}^{[2]} @. \lambda_{\{2\}}^{[2]}                         \\
\frac{\Gamma^{\ptr 2}}{2!} @. \lambda_{\{1\}}^{[1]} @. \lambda_{\{2\}}^{[1]} @. \lambda_{\{1,\, 1\}} \\
\Gamma                     @. \lambda_{\{1\}}     @. \lambda_{\{2\}}                             \\
\{\,o\,\}
\endCD } $ \vskip-3mm     

\vskip-10mm \rightline{picture (4)} \vskip6mm

\noindent
where, in this case (g=6), we have
$$
\lambda_{\{1,\, 1\}} \quad = \quad \CC_{(1)} \ptr \CC_{(1)} + 3 \CC_{(2)} \ptr \CC_{(0)} \ .
\tag{\lma}
$$
The possibilities that may occur are the following (in the sense that they are algebraically admissible
as explained in \S5):

\bitem
(a) \quad $ \CC_{(1)} = 0                  \qquad (\Rightarrow \ \CC_{(2)} = \lambda_{\{1,\,1\}} = 0 \, , \ \dim \RR(\CC) = 7) $

(b) \quad $ \CC_{(1)} \ne 0 \, , \ \CC_{(2)} = \lambda_{\{1,\,1\}} = 0               \qquad (\Rightarrow \ \dim \RR(\CC) = 11) $

(c) \quad $ \CC_{(1)} \ne 0 \, , \ \CC_{(2)} = 0 \, , \ \lambda_{\{1,\,1\}} \ne 0    \qquad (\Rightarrow \ \dim \RR(\CC) = 12) $

(d) \quad $ \CC_{(1)} \ne 0 \, , \ \CC_{(2)} \ne 0 \, , \ \lambda_{\{1,\,1\}} = 0    \qquad (\Rightarrow \ \dim \RR(\CC) = 14) $

(e) \quad $ \CC_{(1)} \ne 0 \, , \ \CC_{(2)} \ne 0 \, , \ \lambda_{\{1,\,1\}} \ne 0  \qquad (\Rightarrow \ \dim \RR(\CC) = 15) $
\eitem

Again, in view of remark (\rbe), we are left to write down basic Pontryagin products.
By (\rijb), non trivial ones
(other then those one of the type $ \, \lambda_{\emptyset} \ptr \lambda_I \equiv \lambda_I \, ) \ $ are

\bitem
$ \lambda_{\{1\}} \ptr \lambda_{\{1\}} \, = \, \lambda_{\{1,1\}} $ \quad in case (c);

$ \lambda_{\{1\}} \ptr \lambda_{\{1\}} \, = \, -3 \lambda_{\{2\}}^{[1]} $ \quad in case (d);

$ \lambda_{\{1\}} \ptr \lambda_{\{1\}} \, = \, -3 \lambda_{\{2\}}^{[1]} + \lambda_{\{1,1\}} $ \quad in case (e).
\eitem

We may note that in this case the list of all products is not so long so we can be more explicit:
beside the products we wrote all the non-trivial products are
$ \ \lambda_{\emptyset}^{[m]} \ptr \lambda_I^{[h]} = {m+h \choose m} \lambda_I^{[m+h]} \, , \ $ for all cases;
$ \ \lambda_{\{1\}}^{[1]} \ptr \lambda_{\{1\}} \, = \, -6 \lambda_{\{2\}}^{[2]} \, , \ $ in cases (d) and (e).
\ Furthermore,
$ \ \lambda_{\emptyset}^{[6-m]} \nts \lambda_I^{[6-h]} = {m+h \choose m} \lambda_I^{[6-m-h]} \, , \ $
$ \ \lambda_{\{1\}}^{[3]} \ptr \lambda_{\{1\}}^{[3]} \, = \, 3 \lambda_{\{2\}}^{[1]} + \lambda_{\{1,1\}} \, , \ $
$ \ \lambda_{\{1\}}^{[3]} \ptr \lambda_{\{1\}}^{[2]} \, = \, 6 \lambda_{\{2\}}^{[2]} \, ; $
\ all other intersection products are trivial.

\newpage

Case $ \, g = 7 \, . \ $

Again by property (\csp), we have $ \, \CC_{(s)} \, = 0 \, $ for $ \, s \ge 4 \, . \ $
Thus, a priori, the basic cycles we should consider, namely those $ \, \lambda_I $'s satisfying $ \ g-s(I)-2d(I) \, \ge \, 0 \, , \ $ are
$ \ \lambda_{\{\emptyset\}},\, \lambda_{\{1\}},\, \lambda_{\{2\}},\, \lambda_{\{3\}},\, \lambda_{\{1,\,1\}},\, \lambda_{\{2,\,1\}} \, . \ $
As usual, the matrices representing the Pontryagin and the intersection products are governed by basic Pontryagin products,
(\mpt), (\mpg) and (\mpb).
As for basic Pontryagin products, in view of formulae (\rijb), we must have
$ \ \lambda_{\emptyset} \ptr \lambda_I \, \equiv \, \lambda_I $ \ and the following relations:

\bitem
$ \ \lambda_{\{1\}} \ptr \lambda_{\{1\}} \, = \, -2 \lambda_{\{2\}}^{[1]} + \lambda_{\{1,1\}} \ ; $

$ \ \lambda_{\{1\}} \ptr \lambda_{\{2\}} \, = \, - 5 \lambda_{\{3\}}^{[1]} + \lambda_{\{2,1\}} \ ; $

$ \ \lambda_{\{1\}} \ptr \lambda_{\{1,1\}} \, = \, a \lambda_{\{3\}}^{[2]} \ , $ for some $ \ a \, \in \, \Q \ $
(in case $ \, \CC_{(3)} \ne 0). $
\eitem

\noindent
(all other basic products are trivial, in fact they are $ \, d $-dimensional cycles of Beauville's degree $ \, s \, $
satisfying $ \ 2d + s \, > \, 7 \, ). $

\rmk{Claim} $ \ \lambda_{\{2,1\}} = 0 \ ; \quad a = 70 \ $ (in case $ \, \CC_{(3)} \ne 0 ) \ ; \quad
\lambda_{\{1,1\}} = 0 \ $ implies $ \ \lambda_{\{3\}} = 0 \, . $

\proof{}
By remark (\alsc), we must have
$$
\big(\lambda_{\{1\}} \ptr \lambda_{\{1\}} \ptr \lambda_{\{1\}}\big)\nts\Theta \quad = \quad
- 18 \lambda_{\{2,\,1\}} \ + \ 90 \lambda_{\{3\}}^{[1]} \ .
$$
On the other hand, using basic products above, (\mpt), (\mpg) and (\mpb) it is easy to check that
$$
\big(\lambda_{\{1\}} \ptr \lambda_{\{1\}} \ptr \lambda_{\{1\}}\big)\nts\Theta \quad = \quad (a + 20) \lambda_{\{3\}}^{[1]} \ .
$$
This gives $ \, \lambda_{\{2,\,1\}} \, = \, (a-70) \lambda_{\{3\}}^{[1]} \, = \, 0 \, . \ $
For $ \, \lambda_{\{3\}} \ne 0 \, , \ $
we obtain $ \, a = 70 \, . \ $ Thus,
because of the relation $ \, \lambda_{\{1\}} \ptr \lambda_{\{1,1\}} \, = \, 70 \lambda_{\{3\}}^{[2]} \, \ne \, 0 \, , \ $
it is also clear that we must have $ \ \lambda_{\{1,1\}} \ne 0  \, . $
\qed

Thus, the picture is the following         \vskip41mm \rightline{picture (5)} \vskip-45mm

\vskip-3mm$ \boxed{ \CD   
\hskip6mm J(\CC) \hskip6mm @. \hskip20mm          @. \hskip20mm          @. \hskip20mm          @. \hskip20mm                \\
\Theta                                                                                                                       \\
\frac{\Gamma^{\ptr 5}}{5!} @. \lambda_{\{1\}}^{[4]}                                                                            \\
\frac{\Gamma^{\ptr 4}}{4!} @. \lambda_{\{1\}}^{[3]} @. \lambda_{\{2\}}^{[3]}                                                     \\
\frac{\Gamma^{\ptr 3}}{3!} @. \lambda_{\{1\}}^{[2]} @. \lambda_{\{2\}}^{[2]} @. \lambda_{\{3\}}^{[2]} @. \lambda_{\{1,\, 1\}}^{[1]}  \\
\frac{\Gamma^{\ptr 2}}{2!} @. \lambda_{\{1\}}^{[1]} @. \lambda_{\{2\}}^{[1]} @. \lambda_{\{3\}}^{[1]} @. \lambda_{\{1,\, 1\}}      \\
\Gamma                     @. \lambda_{\{1\}}     @. \lambda_{\{2\}}     @. \lambda_{\{3\}}                                  \\
\{\,o\,\}
\endCD } $ \vskip-3mm     

\noindent
where, in this case (g=7), by formulae (\rij) we have
$$
\eqalign{
& \lambda_{\{1,\, 1\}} \quad = \quad \CC_{(1)} \ptr \CC_{(1)} \ + \ 2 \, \CC_{(2)} \ptr \CC_{(0)} \ ; \cr
& \lambda_{\{2,\, 1\}} \quad = \quad \CC_{(2)} \ptr \CC_{(1)} \ + \ 5 \, \CC_{(3)} \ptr \CC_{(0)} \ .
}
\tag{\lmb}
$$

As the vanishing of $ \, \lambda_{\{1,1\}} $ \, implies the vanishing of $ \, \lambda_{\{3\}} \, , \ $
also in view of remark (\alsb), it is easy to check that the possibilities that may occur are the following:

\bitem
(a) \quad $ \CC_{(1)} = 0
\qquad (\Rightarrow \ \CC_{(2)} = \CC_{(3)} = \lambda_{\{1,\,1\}} = 0 \, , \ \dim \RR(\CC) = 8) $

(b) \quad $ \CC_{(1)} \ne 0 \, , \ \CC_{(2)} = \lambda_{\{1,\,1\}} = 0
\qquad (\Rightarrow \ \CC_{(3)} = 0 \, , \ \dim \RR(\CC) = 13) $

(c) \quad $ \CC_{(1)} \ne 0 \, , \ \lambda_{\{1,\,1\}} \ne 0 \, , \ \CC_{(2)} = 0
\qquad (\Rightarrow \ \CC_{(3)} = 0 \, , \ \dim \RR(\CC) = 15) $

(d) \quad $ \CC_{(1)} \ne 0 \, , \ \CC_{(2)} \ne 0 \, , \ \CC_{(3)} = \lambda_{\{1,\,1\}} = 0
\qquad (\Rightarrow \ \dim \RR(\CC) = 17) $

(e) \quad $ \CC_{(1)} \ne 0 \, , \ \CC_{(2)} \ne 0 \, , \ \lambda_{\{1,\,1\}} \ne 0 \, , \ \CC_{(3)} = 0
\qquad (\Rightarrow \ \dim \RR(\CC) = 19) $

(f) \quad $ \CC_{(1)} \ne 0 \, , \ \CC_{(2)} \ne 0 \, , \ \CC_{(3)} \ne 0 \, , \ \lambda_{\{1,\,1\}} \ne 0
\qquad (\Rightarrow \ \dim \RR(\CC) = 22) $
\eitem

\noindent
In particular, the admissible dimensions are 8, 13, 15, 17, 19, 22.
Furthermore, all the cases above are algebraically admissible in the sense of section \S 5.

\vskip4mm

Case $ \, g = 8 \, . \ $

Again by property (\csp), we have $ \, \CC_{(s)} \, = 0 \, $ for $ \, s \ge 4 \, . \ $
Thus, a priori, the basic cycles we should consider, namely those $ \, \lambda_I $'s satisfying $ \ g-s(I)-2d(I) \, \ge \, 0 \, , \ $ are
$ \ \lambda_{\{\emptyset\}},\, \lambda_{\{1\}},\, \lambda_{\{2\}},\, \lambda_{\{3\}},\,
\lambda_{\{1,\,1\}},\, \lambda_{\{2,\,1\}},\, \lambda_{\{3,\,1\}},\, \lambda_{\{2,\,2\}} \, . \ $
As usual, the matrices representing the Pontryagin and the intersection products are governed by basic Pontryagin products,
(\mpt), (\mpg) and (\mpb). As for basic Pontryagin products, we must have
$ \ \lambda_{\emptyset} \ptr \lambda_I \, \equiv \, \lambda_I \ $ and the following relations:

\bitem
{\bf a) \quad} $ \lambda_{\{1\}} \ptr \lambda_{\{1\}} \ = \ -{3\over 2} \lambda_{\{2\}}^{[1]} + \lambda_{\{1,1\}} \ ; $

{\bf b) \quad} $ \lambda_{\{1\}} \ptr \lambda_{\{2\}} \ = \ -{10\over 3} \lambda_{\{3\}}^{[1]} + \lambda_{\{2,1\}} \ ; $

{\bf c) \quad} $ \lambda_{\{1\}} \ptr \lambda_{\{3\}} \ = \ \lambda_{\{3,1\}} \ ; $

{\bf d) \quad} $ \lambda_{\{2\}} \ptr \lambda_{\{2\}} \ = \ \lambda_{\{2,2\}} \ ; $

{\bf e) \quad} $ \lambda_{\{1\}} \ptr \lambda_{\{1,1\}} \ = \ 20 \lambda_{\{3\}}^{[2]} - {33\over 2}\lambda_{\{2,1\}}^{[1]} \ ; $

{\bf f) \quad} all other basic products are trivial.

\eitem

\proof{}
Properties a), b), c), d) follow by formulae (\rijb).
Also f) is trivial: those basic products we did not write explicitly
must be linear combinations of cycles of type $ \, \lambda_I^{[m]} \, $ satisfying
$ \ 2d(I)+s(I)+m \, > \, 8 \ $ (which are trivial cycles).

As for proving e) we proceed as follows. First, for reason of dimension and degree we can write
$ \, \lambda_{\{1\}} \ptr \lambda_{\{1,1\}} \, = \, a \lambda_{\{3\}}^{[2]} + b \lambda_{\{2,1\}}^{[1]} \, , \ $
for some $ \, a , \, b \, \in \, \Q \, . \ $
By remark (\alsc), we must have
$$
\big(\lambda_{\{1\}} \ptr \lambda_{\{1\}} \ptr \lambda_{\{1\}}\big)\nts\Theta \quad = \quad
- 18 \lambda_{\{2,\,1\}} \ + \ 60 \lambda_{\{3\}}^{[1]}
$$
On the other hand, using basic products above, (\mpt), (\mpg) and (\mpb) it is easy to check that
$$
\big(\lambda_{\{1\}} \ptr \lambda_{\{1\}} \ptr \lambda_{\{1\}}\big)\nts\Theta \quad = \quad
(b-{\ssize{3\over 2}}) \lambda_{\{2,\,1\}} \ + (2a + 20) \lambda_{\{3\}}^{[1]} \ .
$$
Then we are done by comparing these relations.
\qed

\rmk{Claim}
$ \ \lambda_{\{2,2\}} \ = \ -{10 \over 3} \lambda_{\{3,1\}} \ . $

\proof{}
By remark (\alsc), we must have
$$
\big(\lambda_{\{2\}} \ptr \lambda_{\{1\}} \ptr \lambda_{\{1\}}\big)\nts\Theta \quad = \quad
- 6 \, \lambda_{\{2,\,2\}} \ - 20 \, \lambda_{\{3,\,1\}} \ .
$$
On the other hand, as there are no non-trivial cycles of the form $ \ \sum \lambda_I^{[3-d(I)]} \ $ with $ \ s(I) = 4 \, , \ $
we must have $ \ \lambda_{\{2\}} \ptr \lambda_{\{1\}} \ptr \lambda_{\{1\}} \ = \ 0 \, . \ $
Then we are done.
\qed

Listing the possibilities that may occur, it is useful to observe once for all that
the results above about basic products have a few consequences:
the vanishing of $ \, \lambda_{\{1,1\}} $ \, implies the vanishing of both
$ \, \lambda_{\{3\}} \, $ and $ \, \lambda_{\{2,1\}} \, ; \ $
the non-vanishing of $ \, \lambda_{\{3,1\}} \, $ (equivalently of $ \, \lambda_{\{2,2\}} \, ) \ $
implies the non-vanishing of $ \, \lambda_{\{3\}} \, $ and $ \, \lambda_{\{1,1\}} \, . \ $
Eventually, it is easy to check that the possibilities that may occur are the following:

\bitem
(a) \quad $ \CC_{(1)} = 0
\qquad (\Rightarrow \ \CC_{(2)} = \CC_{(3)} = \lambda_{\{1,\,1\}} = \lambda_{\{2,\,1\}} = 0 \, , \ \dim \RR(\CC) = 9) $

(b) \quad $ \CC_{(1)} \ne 0 \, , \ \CC_{(2)} = \lambda_{\{1,\,1\}} = 0
\qquad (\Rightarrow \ \CC_{(3)} = \lambda_{\{2,\,1\}} = 0 \, , \ \dim \RR(\CC) = 15) $

(c) \quad $ \CC_{(1)} \ne 0 \, , \ \lambda_{\{1,\,1\}} \ne 0 \, , \ \CC_{(2)} = 0
\qquad (\Rightarrow \ \CC_{(3)} = \lambda_{\{2,\,1\}} = 0 \, , \ \dim \RR(\CC) = 18) $

(d) \quad $ \CC_{(1)} \ne 0 \, , \ \CC_{(2)} \ne 0 \, , \ \CC_{(3)} = \lambda_{\{1,\,1\}} = 0
\qquad (\Rightarrow \ \lambda_{\{2,\,1\}} = 0 \, , \ \dim \RR(\CC) = 20) $

(e) \quad $ \CC_{(1)} \ne 0 \, , \ \CC_{(2)} \ne 0 \, , \ \lambda_{\{1,\,1\}} \ne 0 \, , \ \CC_{(3)} = \lambda_{\{2,\,1\}} = 0
\qquad (\Rightarrow \ \dim \RR(\CC) = 23) $

(f) \quad $ \CC_{(1)} \ne 0 \, , \ \CC_{(2)} \ne 0 \, , \ \lambda_{\{1,\,1\}} \ne 0 \, , \ \lambda_{\{2,\,1\}} \ne 0 \, , \ \CC_{(3)} = 0
\qquad (\Rightarrow \ \dim \RR(\CC) = 25) $

(g) \quad $ \CC_{(1)} \ne 0 \, , \ \CC_{(2)} \ne 0 \, , \ \CC_{(3)} \ne 0 \, , \ \lambda_{\{1,\,1\}} \ne 0 \, , \ \lambda_{\{2,\,1\}} = 0
\qquad (\Rightarrow \ \dim \RR(\CC) = 27) $

(h) \quad $ \CC_{(1)} \ne 0 \, , \ \CC_{(2)} \ne 0 \, , \ \CC_{(3)} \ne 0 \, , \ \lambda_{\{1,\,1\}} \ne 0 \, , \ \lambda_{\{2,\,1\}} \ne 0
\qquad (\Rightarrow \ \dim \RR(\CC) = 29) $

(i) \quad $ \CC_{(1)} \ne 0 \, , \ \CC_{(2)} \ne 0 \, , \ \CC_{(3)} \ne 0 \, , \ \lambda_{\{1,\,1\}} \ne 0 \, , \ \lambda_{\{2,\,1\}} \ne 0
\, , \ \lambda_{\{3,\,1\}} \ne 0 \qquad (\Rightarrow \ \dim \RR(\CC) = 30) $
\eitem

\noindent
In particular, the admissible dimensions are 9, 15, 18, 20, 23, 25, 27, 29, 30 and the picture is the following

\vskip-3mm$ \boxed{ \CD   
\hskip6mm J(\CC) \hskip6mm @. \hskip20mm          @. \hskip20mm          @. \hskip20mm          @. \hskip20mm               @. \hskip20mm  @. \hskip20mm               \\
\Theta                                                                                                                                                  \\
\frac{\Gamma^{\ptr 6}}{6!} @. \lambda_{\{1\}}^{[5]}                                                                                                       \\
\frac{\Gamma^{\ptr 5}}{5!} @. \lambda_{\{1\}}^{[4]} @. \lambda_{\{2\}}^{[4]}                                                                                \\
\frac{\Gamma^{\ptr 4}}{4!} @. \lambda_{\{1\}}^{[3]} @. \lambda_{\{2\}}^{[3]} @. \lambda_{\{3\}}^{[3]} @. \lambda_{\{1,\, 1\}}^{[2]}                             \\
\frac{\Gamma^{\ptr 3}}{3!} @. \lambda_{\{1\}}^{[2]} @. \lambda_{\{2\}}^{[2]} @. \lambda_{\{3\}}^{[2]} @. \lambda_{\{1,\, 1\}}^{[1]} @. \lambda_{\{2,\, 1\}}^{[1]} \\
\frac{\Gamma^{\ptr 2}}{2!} @. \lambda_{\{1\}}^{[1]} @. \lambda_{\{2\}}^{[1]} @. \lambda_{\{3\}}^{[1]} @. \lambda_{\{1,\, 1\}}     @. \lambda_{\{2,\, 1\}}    @. \lambda_{\{3,\, 1\}}     \\
\Gamma                     @. \lambda_{\{1\}}     @. \lambda_{\{2\}}     @. \lambda_{\{3\}}                                                             \\
\{\,o\,\}
\endCD } $ \vskip-3mm     

\vskip5mm \rightline{picture (6)}

\noindent
where, in this case (g=8), in view of formulae (\rij) we have
$$
\eqalign{
& \lambda_{\{1,\, 1\}} \quad = \quad \CC_{(1)} \ptr \CC_{(1)} \ + \ {3\over 2} \, \CC_{(2)} \ptr \CC_{(0)} \ ; \cr
& \lambda_{\{2,\, 1\}} \quad = \quad \CC_{(2)} \ptr \CC_{(1)} \ + \ {10\over 3} \, \CC_{(3)} \ptr \CC_{(0)} \ ; \cr
& \lambda_{\{3,\, 1\}} \quad = \quad \CC_{(3)} \ptr \CC_{(1)} \quad = \quad - {3\over 10} \CC_{(2)} \ptr \CC_{(2)} \ .
}
\tag{\lmc}
$$

\vskip4mm

$ \RR(\CC) \ $ for trigonal curves.

The tautological ring of a trigonal curve was already discussed by Beauville in [Be3].
In terms of our models the situation is very simple:
as $ \, \CC_{(s)} = 0 \, $ for $ \, s \ge 2 \, , \ $ there exist an integer $ \, k \, $ such that the basic
non-trivial cycles are the $ \, \lambda_I $'s with $ \, I \, = \, {\{1, \, ..., \, 1\}} $ and $ \ d(I) \le k \, . $
We may note that as $ \, g - s(I) - 2 d(I) \ge 0 \, $ we have $ \, k \le {g \over 3} \, . \ $
It is also clear that basic Pontryagin products are
$ \, \lambda_I \ptr \lambda_J = \lambda_{I\cup J} \, , \ I $ and $ J $ as above.
In particular, \ $ \dim \RR(\CC) \, = \, (g+1) + (g-2) + (g-5) + ... + (g+1-3k) \ = \ (k+1) \left(g + 1 - {3 \over 2}k\right) \, . $

\vskip4mm

$ \RR(\CC) \ $ for curves with a $ \, g^1_4 \, . $

In this case, $ \, \CC_{(s)} = 0 \, $ for $ \, s \ge 3 \, . \ $ Thus, there exist integers
$ \, k_0 \ge k_1 \ge ...  \ge k_m \, $ such that the basic non-trivial cycles are the
$ \, \lambda_I $'s with $ \, I \, = \, {\{2, \, ..., \, 2, \, 1, \, ...,\, 1 \}} \, , \ $
where \lq\lq 2" is repeated $ \, h $-times $(h\le m)$ and \lq\lq 1" is repeated up to $ \, k_h $-times.
As $ \ s(I) = 2h+k_h \ $ and $ \ d(I) = h+k_h \, , \ $ admissibility of multi indexes reduces to the conditions
$$
\qquad 4 \, h \, + \, 3 \, k_h \ \ \le \ \ g \ , \qquad 0 \ \le \ h \ \le \ m \, ,
$$
i.e $ \ 3 k_0 \le g \, , \ 4 + 3 k_1 \le g \, , \ 8 + 3 k_2 \le g ... \ $ etcetera.
As $ \, s(I) \, $ and $ \, d(I) \, $ determine $ \, I \, , $
\ indeed \lq\lq 2" is repeated $ \, s(I) -d(I) \, $ times and \lq\lq 1" is repeated $ \, 2d(I) -s(I) \, $ times,
also in view of corollary (\pzc), the spaces
$ \, \left\langle \, ..., \, \lambda_I^{[t]} , \, ... \, \right\rangle_{t=0}^{g-s(I)-2d(I)} \, $ are independent, therefore
$$
dim \RR(\CC) \quad = \quad \sum_{h=0}^m \ g - 4 h - 3 k_h
$$


\vskip12mm

\centerline{\bf REFERENCES}

\vskip3mm
\widestnumber\key{MMMM}

\ref
\key ACGH
\by E. Arbarello, M. Cornalba, P.A. Griffith, J. Harris
\book Geometry of algebraic curves
\publ Springer (267)
\vol I
\yr 1985
\endref

\ref
\key Be1
\by A. Beauville
\paper Quelques remarques sur la tranformation de Fourier dans l'anneau de Chow d'une variet\'e ab\'elienne
\jour L.N.M.
\vol 1016
\yr 1983
\pages 238-260
\endref

\ref
\key Be2
\by A. Beauville
\paper Sur l'anneau de Chow d'une variet\'e ab\'elienne
\jour Math. Ann.
\vol 273
\yr 1986
\pages 647-651
\endref

\ref
\key Be3
\by A. Beauville
\paper Algebraic cycles on Jacobian varieties
\jour arXiv:math.AG/0204188v1
\vol
\yr 2002
\endref

\ref
\key Ce
\by G. Ceresa
\paper $ C $ is not algebraically equivalent to $ C^- $ in its Jacobian
\jour Ann. of Math.
\vol 117
\yr 1983
\pages 285-291
\endref

\ref
\key CG
\by E. Colombo, B. van Geemen
\paper Note on curves in a Jacobian
\jour Compositio Math.
\vol 88
\yr 1993
\pages 333-353
\endref

\ref
\key DM
\by C. Deninger, J. Murre
\paper Motivic decomposition of abelian schemes and the Fourier transform
\jour J. Reine Angew. Math.
\vol 422
\yr 1991
\pages 201-219
\endref

\ref
\key Fu
\by W. Fulton
\book Intersection theory
\publ Ergeb. der Math. 3 Folge, Band 2.
\yr 1984, Springer Verlag, Berlin
\endref

\ref
\key GMV
\by M. Green, J. Murre, C. Voisin
\paper Algebraic Cycles and Hodge Theory
\jour Lecture Notes in Mathematics
\vol 1594
\yr 1993
\pages
\endref

\ref
\key K
\by K. K\"unnemann
\paper On the Chow motive of an abelian scheme
\jour Proceedings of the conference on Motives (AMS Proc. in pure Math.)
\vol 88
\yr 1991
\pages 333-353
\endref

\ref
\key Ma
\by G. Marini
\paper Algebraic cycles on abelian varieties and their decomposition
\jour {\bf Bollettino U.M.I.}
\vol (8) 7-B
\yr 2004
\pages 231-240 \
\endref

\ref
\key Muk
\by S. Mukai
\paper Duality between $ D(X) $ and $ D(\hat X) $ with its applications to Picard sheaves
\jour Nagoya Math. J.
\vol 81
\yr 1981
\pages 153-175
\endref

\enddocument

\newpage

\vskip7.5cm

\ bed\, = \, 1\ \   \ fut\, = \, 2\ \   \ tgf\, = \, 3\ \   \ ntt\, = \, 4\ \   \ fte\, = \, 5\ \   \ okt\, = \, 6\ \   \ ppz\, = \, 7\ \   \ csp\, = \, 8\ \
\ ppt\, = \, 9\ \   \ pzn\, = \, 10\ \   \ pzl\, = \, 11\ \   \ pzm\, = \, 12\ \   \ pzc\, = \, 13\ \   \ pzr\, = \, 14\ \   \ pzd\, = \, 15\ \   \ pzv\, = \, 16\ \
\ rbp\, = \, 17\ \   \ mpt\, = \, 18\ \   \ mpg\, = \, 19\ \   \ mpb\, = \, 20\ \   \ cij\, = \, 21\ \   \ rij\, = \, 22\ \   \ rijb\, = \, 23\ \
\ cvn\, = \, 24\ \   \ chij\, = \, 25\ \   \ ami\, = \, 26\ \   \ als\, = \, 27\ \   \ alsb\, = \, 28\ \   \ alsc\, = \, 29\ \    \ rbe\, = \, 30\ \
\ lma\, = \, 31\ \   \ lmb\, = \, 32\ \   \ lmc\, = \, 33\ \   \ alsd\, = \, 34\ \
\ lcc\, = \, 41\ \   \ lcu\, = \, 42\ \   \ ppc\, = \, 43\ \   \ pzz\, = \, 44\ \

\enddocument